\RequirePackage{xstring}
\StrBehind*{\jobname}{.}[\DocOpts]

\wlog{DocOpts: \DocOpts}

\documentclass[\DocOpts]{amsart}

\usepackage{amssymb}
\usepackage{enumerate}
\input xy
\xyoption{all}
\usepackage{ifdraft}
\ifdraft{
\usepackage[breaklinks,bookmarks,final]{hyperref}
}{
}
\providecommand{\pagewiselinenumbers}{\relax}

\newtheorem{theorem}{Theorem}[section]
\newtheorem{lemma}[theorem]{Lemma}
\newtheorem{corollary}[theorem]{Corollary}

\newtheorem{fact}[theorem]{Fact}
\newtheorem{proposition}[theorem]{Proposition}
\newtheorem{propdef}[theorem]{Proposition-Definition}
\newtheorem{claim}[theorem]{Claim}
\theoremstyle{definition} 
\newtheorem{remark}[theorem]{Remark}
\newtheorem{definition}[theorem]{Definition}
\newtheorem{question}[theorem]{Question}

\def\characteristic{\operatorname{char}}
\def\fr{\operatorname{Fr}}
\def\fix{\operatorname{Fix}}

\def\acfa{\operatorname{ACFA}}
\def\acf{\operatorname{ACF}}

\def\cb{\operatorname{Cb}}

\def\tp{\operatorname{tp}}
\def\qftp{\operatorname{qftp}}
\def\qf{\operatorname{qf}}

\def\acl{\operatorname{acl}}
\def\dcl{\operatorname{dcl}}

\def\alg{\operatorname{alg}}
\def\perf{\operatorname{perf}}

\def\aut{\operatorname{Aut}}
\def\bir{\operatorname{Bir}}
\def\loc{\operatorname{loc}}
\def\inv{\operatorname{inv}}
\def\trdeg{\operatorname{trdeg}}

\def\lring{\mathcal L_{\operatorname{ring}}}
\def\lsigma{\mathcal L_\sigma}

\def\C{\mathcal{C}}
\def\G{\mathcal{G}}

\def\U{\mathcal{U}}
\def\I{\mathcal{I}}
\renewcommand{\O}{\mathcal{O}}
\def\Q{\mathcal{Q}}

\def\PP{{\mathbb{P}}}

\def\VV{\mathbb{V}}

\newcommand{\NN}{\mathbb{N}}

\renewcommand{\AA}{\mathbb{A}}
\newcommand{\LL}{\mathbb{L}}

\newcommand{\isom}{\cong}

\providecommand{\id}{\operatorname{id}}

\newcommand{\ra}[1][]{\xrightarrow{#1}}

\newcommand{\dom}[1]{{\operatorname{dom} (#1)}}

\newcommand{\prolong}[2][]{\nabla_{#1}(#2)}

\newcommand{\dto}{\dashrightarrow}

\def\Ind#1#2{#1\setbox0=\hbox{\(#1x\)}\kern\wd0\hbox to 
0pt{\hss\(#1\mid\)\hss}\lower.9\ht0\hbox to 
0pt{\hss\(#1\smile\)\hss}\kern\wd0}
\def\ind{\mathop{\mathpalette\Ind{}}}
\def\notind#1#2{#1\setbox0=\hbox{\(#1x\)}\kern\wd0\hbox to 
0pt{\mathchardef\nn=12854\hss\(#1\nn\)\kern1.4\wd0\hss}\hbox to 
0pt{\hss\(#1\mid\)\hss}\lower.9\ht0 \hbox to 
0pt{\hss\(#1\smile\)\hss}\kern\wd0}
\def\nind{\mathop{\mathpalette\notind{}}}

\newenvironment{claimproof}[1]{%
\par\vspace{.1in}\noindent\emph{Proof of Claim:}\space#1}{%
\vspace{.1in}\hfill\(\maltese\)}

\date{\today}

\title{Binding groups for algebraic dynamics}
\author{Moshe Kamensky}
\address{Moshe Kamensky\\
Ben-Gurion University of the Negev\\
Department of Mathematics\\
Be'er-Sheva, 8410501\\
Israel}
\email{kamensky.bgu@gmail.com}
\author{Rahim Moosa}
\address{Rahim Moosa\\
University of Waterloo\\
Department of Pure Mathematics\\
200 University Avenue West\\
Waterloo, Ontario \  N2L 3G1\\
Canada}
\email{rmoosa@uwaterloo.ca}

\thanks{R. Moosa was partially supported by an NSERC Discovery Grant.}

\begin{document}


\begin{abstract}
	A binding group theorem is proved in the context of quantifier-free 
	internality to the fixed field in the theory \(\acfa_0\).
	This is articulated as a statement about the birational geometry of 
	isotrivial algebraic dynamical systems, and more generally isotrivial 
	\(\sigma\)-varieties.
	It asserts that if \((V,\phi)\) is an isotrivial \(\sigma\)-variety then a 
	certain subgroup of the group of birational transformations of~\(V\), 
	namely those that preserve all the relations between \((V,\phi)\) and the 
	trivial dynamics \((\AA^1,\id)\), is in fact an algebraic group.
	Several application are given including new special cases of the Zariski 
	Dense Orbit Conjecture and the Dixmier-Moeglin Equivalence Problem in 
	algebraic dynamics, as well as finiteness results about the existence of 
	nonconstant invariant rational functions on cartesian powers of 
	\(\sigma\)-varieties. These applications give algebraic-dynamical analogues 
	of recent results in differential-algebraic geometry.
\end{abstract}

\maketitle

\setcounter{tocdepth}{1}
\tableofcontents

\section{Introduction}

\pagewiselinenumbers%

\noindent
Fix an algebraically closed field~\(k\) of characteristic zero, equipped with an automorphism~$\sigma$.
This paper is concerned with the  birational algebraic geometry, and model 
theory, of \emph{rational \(\sigma\)-varieties} over \((k,\sigma)\); that is, irreducible algebraic varieties~\(V\) over~\(k\) equipped with a dominant 
rational map \(\phi:V\dto V^\sigma\).
Note that such a structure on~\(V\) determines, and is determined by, an 
endomorphism of the rational function field~\(k(V)\) that extends~\(\sigma\); namely, given $f\in k(V)$, precompose $f^\sigma\in k(V^\sigma)$ with $\phi$ to obtain $f^\sigma\circ\phi\in k(V)$.

The algebraic dynamics literature usually only considers the autonomous 
situation where~\(\sigma\) is the identity on~\(k\), and hence 
\(V=V^\sigma\). These are \emph{rational dynamical systems.}
While they are an important special case for us too, we work generally in the 
possibly nonautonomous context.\footnote{One reason to not restrict attention 
only to the autonomous case is to allow the taking of generic fibres: even if 
both \((V,\phi)\) and \((W,\psi)\) are rational dynamical systems, the 
generic fibre of an equivariant map \(g:(V,\phi)\to(W,\psi)\), namely the 
induced base extension of~\(V\) to~\(k(W)\), should be considered with its 
natural nonautonomous \(\sigma\)-variety structure coming from the 
endomorphism of \(k(W)\) induced by~\(\psi\).}

The natural subobjects here are the \emph{invariant} subvarieties; namely, 
irreducible subvarieties \(X\subseteq{}V\) over~\(k\) on which \(\phi\) restricts to a 
dominant rational map from \(X\) to \(X^\sigma\).
And the morphisms of this category,  \(g:(V_1,\phi_1)\dto (V_2,\phi_2)\), are 
the dominant rational maps \(g:V_1\to V_2\) that are \emph{equivariant} in 
the sense that \(\phi_2\circ g=g^\sigma\circ\phi_1\).

Our focus is on \emph{isotrivial} 
\(\sigma\)-varieties: those that are, after base extension, equivariantly birationally equivalent to 
a trivial \(\sigma\)-variety, namely a variety equipped with the identity transformation.
More precisely, a rational \(\sigma\)-variety \((V,\phi)\) is isotrivial
if there is a commuting diagram of the form
\[
\xymatrix{
(V\times Z,\phi\times\psi)\ar[dr]\ar@{-->}[rr]^{g}&& 
(\AA^\ell\times Z,\id\times\psi)\ar[dl]\\
&(Z,\psi)
}
\]
where \((Z,\psi)\) is another rational \(\sigma\)-variety over~\(k\) and 
\(g\) is birational onto its image.
A basic example of a nontrivial but isotrivial rational dynamical system is 
the map \(\phi:\AA^1\to\AA^1\) given by \(\phi(x)=x+1\).
The trivialisation is obtained by taking \((Z,\psi)\) to be \((\AA^1,\phi)\) 
itself, and~\(g\) to be given by \(g(x,y)=(x-y,y)\).
In fact, isotrivial \(\sigma\)-varieties are ubiquitous; see, for example, Corollary~\ref{nonortho-qfint-geo} below, which says that if a \(\sigma\)-variety admits any nonvacuous algebraic relation to trivial dynamics, even after base extension, then already without base extension it admits a positive-dimensional isotrivial image.
An even more convincing example of the centrality of isotrivial \(\sigma\)-varieties is the Zilber dichotomy (a deep model-theoretic result established in~\cite{acfa}) which implies, roughly speaking, that as soon as some cartesian power of \((V,\phi)\) admits a sufficiently rich algebraic family of invariant subvarieties, then \((V,\phi)\) admits a positive-dimensional isotrivial image.

Isotrivial rational dynamics were studied (by very different means) 
in~\cite{bms}, where they were shown to always come from the action of an 
algebraic group, very much like in the above example where the relevant group 
action is that of the additive group on the affine line.
We are partly motivated by the desire to extend that work to the possibly 
nonautonomous setting, and to give a model-theoretic account.
Nevertheless, our approach gives significant new information even for 
rational dynamical systems.

Our main result is that a certain natural group of birational transformations 
of an isotrivial \(\sigma\)-variety is in fact an algebraic group.
To describe the result we need some notation.
First of all, for any field extension \(K\supseteq k\), let us denote by 
\(\bir_K(V)\) the group of birational transformations of~\(V\) over~\(K\).
That is, birational maps \(\delta:V_K\dto V_K\) under composition.
We are interested in ``algebraic subgroups" of \(\bir_k(V)\) in the following 
sense:

\begin{definition}\label{def:agbt}
	By an \emph{algebraic group of birational transformations of~\(V\)} we mean 
	an algebraic group~\(G\) over~\(k\) equipped with a rational map 
	\(\theta:G\times V\dto V\) over~\(k\) that determines an injective group 
	homomorphism \(G(K)\to \bir_K(V)\), for any field extension \(K\supseteq 
	k\), given by \(w\mapsto\theta_w\).
\end{definition}

\begin{remark}
\begin{itemize}
\item[(a)]
The algebraic group~$G$ in this definition need not be connected.
As such, one has to be a little more careful in defining rational maps, see~$\S$\ref{sec:ag} below for our algebro-geometric conventions.

\item[(b)]
Note that~\(\theta\) 
	makes~\(V\) into a pre-homogeneous variety for~\(G\), in the sense of 
	Weil~\cite{weil}, and it follows that 
	after replacing~\(V\) by a birationally equivalent copy, we get an honest 
	regular (rather than rational) algebraic group action.
	Actually, as~$G$ may not be connected, we are really using Zaitsev's modern treatment in~\cite{zaitsev} of Weil's regularisation theorem, to which we refer the reader for details.
\end{itemize}
\end{remark}

For example, while the usual action of \(\mathbb G_a\) on \(\AA^1\) is an 
algebraic group of automorphisms, the variant given by 
\((g,x)\mapsto\frac{x}{gx+1}\) is only an algebraic group of birational 
transformations.\footnote{This example is slightly artificial since it comes 
from a regular action on \(\PP^1\), but in higher dimensions there are algebraic groups of birational transformations of a projective variety that are not algebraic groups of automorphisms.}
These are well-studied objects in birational algebraic geometry, especially 
in the case when \(V=\PP^n\); See, for example,~\cite{blanc} and the 
references therein.

The following abstract subgroup of \(\bir(V)\) captures the interaction 
between a given rational \(\sigma\)-variety structure on~\(V\) and trivial 
dynamics.
It appears to us to be a fundamental object in algebraic dynamics that has 
not been studied before, even in the autonomous case.

\begin{definition}[Binding group of a \(\sigma\)-variety]\label{def:bgrds}
	Fix an irreducible rational \(\sigma\)-variety \(\VV:=(V,\phi)\) 
	over \((k,\sigma)\), and let \(\LL:=(\AA^1,\id)\) denote the trivial 
	dynamics on the affine line.
	Let
	\[
	\I_{\VV}:=
	\left\{
	\begin{tabular}{ll}
		irreducible invariant subvarieties of \(\VV^r\times\LL^s\) over~\(k\) 
		that project\\
		dominantly onto each copy of~\(V\), for all \(r\geq 1\) and \(s\geq 0\)\\
	\end{tabular}
	\right\}
	\]
	Fix, now, a field extension \(K\supseteq k\).
	For each \(r\geq 1\) and \(s\geq 0\), embed \(\bir_K(V)\) into 
	\(\bir_K(V^r\times\AA^s)\) by acting diagonally on~\(V^r\) and trivially on 
	\(\AA^s\).
	Let
	\[
	\bir_K(\VV/\LL):=
	\left\{
	\delta\in\bir_K(V):
	\begin{tabular}{ll}
		if \(X\in\I_{\VV}\) then \(\delta\) restricts to a\\
		birational transformation of \(X_K\)\\
	\end{tabular}
	\right\}
	\]
	That is, \(\bir_K(\VV/\LL)\) is the group of birational transformations 
	of~\(V\) over~\(K\) that preserve all invariant \(k\)-definable algebraic 
	relations between cartesian powers of \(\VV\) and~\(\LL\).
\end{definition}

We show that the \emph{a priori} abstract binding group of an isotrivial 
\(\sigma\)-variety is in fact an algebraic group:

\begin{theorem}\label{thm:bgrds}
	Suppose 
	\(\VV:=(V,\phi)\) is a rational \(\sigma\)-variety over~\((k,\sigma)\).
	If \(\VV\) is isotrivial then there exists an algebraic group~\(G
	\) of birational transformations of~\(V\) such that 
	\(\bir_K(\VV/\LL)=G(K)\), for any field extension \(K\supseteq k\).
\end{theorem}

The proof of this theorem appears in Section~\ref{sect:bg} below.

We delay discussion of applications of this theorem to later in the 
Introduction, addressing first its model-theoretic formulation.
The model-theorist will by now have realised that Theorem~\ref{thm:bgrds} has 
something to do with what is often called ``the binding group theorem".
Indeed, what we prove is a quantifier-free binding group theorem for the 
theory of difference-closed fields (\(\acfa\)), introduced by Chatzidakis and 
Hrushovski in~\cite{acfa}.
Recall that a \emph{difference field} is a field equipped with an 
endomorphism, and it is \emph{difference-closed} if every system of algebraic
\emph{difference equations} (that is, polynomial equations in variables 
\(x,\sigma(x),\sigma^2(x),\dots\)) that is consistent -- in the sense that it 
has a solution in some difference field extension --  has a solution.
The connection to algebraic dynamics is that a \(\sigma\)-variety 
\((V,\phi)\) can be seen as encoding the first-order difference equation 
\(\sigma(x)=\phi(x)\).
That is, given any difference field extension 
\((K,\sigma)\supseteq(k,\sigma)\), we can consider those \(K\)-points 
of~\(V\) on which \(\sigma\)~and~\(\phi\) agree.
This is a quantifier-free definable set in \((K,\sigma)\) that we denote by 
\({(V,\phi)}^\sharp(K)\).
The fact that the class of difference-closed fields is axiomatisable (by 
\(\acfa\)) means that, instead of considering all possible difference field 
extensions of \((k,\sigma)\), we can work in a fixed large difference-closed 
extension, \(\U\), that serves as a universal domain for difference-algebraic 
geometry.

We associate to each rational \(\sigma\)-variety \((V,\phi)\) over 
\((k,\sigma)\), the quantifier-free type~\(q(x)\) over~\(k\) which asserts 
that~\(x\) is a Zariski generic point of~\(V\) over~\(k\) and that 
\(\sigma(x)=\phi(x)\).
This turns out to be a complete quantifier-free type which we call the 
\emph{generic type} of \((V,\phi)\).
We call those complete quantifier-free types that arise in this way, 
\emph{rational types}.
The model theory of~\(q\) in \(\U\) will control, and is controlled by, the 
birational geometry of \((V,\phi)\).
In particular, in characteristic zero, \((V,\phi)\) being isotrivial 
corresponds precisely to \(q\) being \emph{quantifier-free internal} to the 
fixed field.
Quantifier-free internality in \(\acfa\) was introduced in~\cite{acfa} and is 
discussed  at length in~\(\S\)\ref{subsect:qfint} below.
The model-theoretic content of this paper has to do, therefore, with the 
structure of rational types in \(\acfa\) that are quantifier-free internal to 
the fixed field.
In particular, we introduce a binding group:

\begin{definition}[Binding group of a rational type]\label{def:bgrt}
	Given a rational quantifier-free type \(q\) over~\((k,\sigma)\), we denote 
	by \(\aut_{\qf}(q/\fix(\sigma))\) the (abstract) subgroup of 
	permutations~\(\delta\) of \(q(\U)\) satisfying:
	\begin{equation*}
		\theta(a,c)\text{ holds }\iff \theta(\delta(a),c)\text{ holds}
	\end{equation*}
	for any (quantifier-free) formula \(\theta(x,y)\) over~\(k\) in the 
	language of rings, any tuple~\(a\) of realisations of~\(q\), and any 
	tuple~\(c\) of elements of the fixed field.
\end{definition}

And we prove a binding group theorem:

\begin{theorem}\label{thm:bgrt}
	Working in a sufficiently saturated difference-closed field~\(\U\) of 
	characteristic zero, suppose \(q\) is a rational quantifier-free type 
	over~\(k\) that is quantifier-free internal to the fixed field.
	There exists a quantifier-free definable group \(\G\) over~\(k\), with a 
	relatively quantifier-free definable action on \(q(\U)\) over~\(k\), such 
	that~\(\G\) and \(\aut_{\qf}(q/\fix(\sigma))\) are isomorphic as groups 
	acting on \(q(\U)\).

	In fact, if \(q\) is the generic type of an isotrivial rational 
	\(\sigma\)-variety \((V,\phi)\) over~\((k,\sigma)\), and \(\theta:G\times 
	V\dto V\) is the algebraic group of birational transformations of~\(V\) 
	given by Theorem~\ref{thm:bgrds}, then there is an isomorphism \(\rho:G\to 
	G^\sigma\) of algebraic groups over~\(k\) such that 
	\(\G={(G,\rho)}^\sharp(\U)\) and \(\theta\) restricts to the action 
	of~\(\G\) on \(q(\U)\).
\end{theorem}

This is also proved in Section~\ref{sect:bg}.

Let us briefly recall the model-theoretic precedents to this theorem.
A crucial aspect of usual internality in totally transcendental theories, of 
a complete type~\(p\), say, to a definable set~\(X\), is that the witness to 
internality may involve more parameters than those over which~\(p\) and~\(X\) 
are defined.
This dependence on additional parameters is controlled by the binding group 
(or \emph{liaison} group), a definable group acting definably on the 
realisations of~\(p\) and agreeing with the action of the group of 
automorphisms of the universe that fix~\(X\) pointwise.
The existence and importance of the binding group was already recognised by 
Zilber~\cite{zilber} in the late nineteen-seventies.
Poizat~\cite{poizat} realised that when applied to differentially closed 
fields, binding groups recover Kolchin's differential Galois theory.
Hrushovski developed the subject in its current form, first working with 
stable theories but eventually in complete generality: in~\cite{udigeneral} 
the binding group is constructed (as a type-definable group) from internality 
assuming only that the set \(X\) is stably embedded.
Based on Hrushovski's construction, the first author, in~\cite{Moshe}, 
extended the theory of binding groups to the quantifier-free fragment (or 
indeed arbitrary fragments) of a theory, very much with \(\acfa\) in mind.
The focus of~\cite{Moshe} is linear difference equations and the development 
of a binding group theory that recovers the difference Galois theory of Van 
der Put and Singer.
Moreover, it is concerned with the internality of one definable set in 
another, and not of generic types.
In particular, the results there do not immediately apply to the birational 
geometry of rational \(\sigma\)-varieties.
Nevertheless, while we do not directly rely on~\cite{Moshe}, that work should 
be considered the immediate predecessor of this one, and it very much 
influences our construction.

\subsection{Applications}
We now describe several applications.
The proofs, and more detailed statements,  of the following theorems appear 
in Section~\ref{sect:applications} below.
Each of these applications has both a formulation in terms of the birational 
geometry of rational \(\sigma\)-varieties, as well as in terms of the model 
theory of rational types in~\(\acfa_0\). Here, in the Introduction, we focus 
on the geometric formulations.

First of all, if we restrict attention to rational dynamics \(\phi:V\dto V\), so when $\sigma$ is the identity on~$k$, then we recover some of the main results of~\cite{bms}.
In particular, we 
are able to show that if \((V,\phi)\) is isotrivial then~\(\phi\) comes from 
an algebraic group action; this is~\cite[Corollary~A]{bms} and appears as 
Theorem~\ref{thm:translational} below. The proof in~\cite{bms} is somewhat 
involved and computational, using mostly elementary methods from algebraic 
dynamics.  We deduce it here by observing that in the autonomous isotrivial 
case, \(\phi\in\bir(\VV/\LL)\), and so the algebraic group is the one given 
by Theorem~\ref{thm:bgrds}.

As a more or less standard corollary, we resolve a special case of the 
\emph{Zariski dense orbit conjecture} (from~\cite{alicetom}) that we don't 
think has been observed before (though it follows from the results of~\cite{bms} as well):

\begin{theorem}[Appearing as Corollary~\ref{cor:translational} below]
	Suppose \(\phi:V\to V\) is an automorphism of an algebraic variety over~\(k\) 
	such that \((V,\phi)\) is isotrivial.
	If \((V,\phi)\) admits no nonconstant invariant rational functions then 
	there is \(a\in V(k)\) such that the orbit of~\(a\) under~\(\phi\) is 
	Zariski dense in~\(V\).
\end{theorem}

We are, similarly, able to give a model-theoretic 
account of~\cite[Corollary~B]{bms} in Corollary~\ref{cor:nonorthbound-auto} 
below, but we leave the formulation of that result for later in the 
Introduction.  Explaining the results of~\cite{bms} from the point of view of 
model-theoretic binding groups was one of the motivations for this work.

Our first new application is about the number of \emph{maximal} proper 
invariant subvarieties. A necessary condition for a rational 
\(\sigma\)-variety, \((V,\phi)\) over~\((k,\sigma)\), to admit only finitely 
many maximal proper invariant subvarieties over~\(k\) is that \((V,\phi\)) 
admit no nonconstant invariant rational functions.
Here, an \emph{invariant rational function} on \((V,\phi)\) is a rational 
function on~\(V\) over~\(k\) that is fixed by the endomorphism of \(k(V)\) 
that~\(\phi\) induces.
Such a rational function would, by taking level sets, give rise to infinitely 
many distinct codimension~\(1\), and hence maximal proper, invariant 
subvarieties.
The question of whether this condition is sufficient -- that is, whether 
having no nonconstant invariant rational functions implies having only 
finitely many maximal proper invariant subvarieties --  is sometimes called 
the \emph{Dixmier-Moeglin equivalence problem} in algebraic dynamics, at 
least in the case when \(\phi\) is an automorphism of a projective variety 
(see~\cite[Conjecture~8.5]{BRS} and also~\cite{dme-survey} for a survey of 
Dixmier-Moeglin-type problems).
Using binding groups, we resolve the problem for isotrivial 
\(\sigma\)-varieties:

\begin{theorem}[Appearing as Theorem~\ref{dme} below]\label{thm:dme}
	Suppose 
	\((V,\phi)\) is an isotrivial rational \(\sigma\)-variety over 
	\((k,\sigma)\).
	If \((V,\phi)\) has no nonconstant invariant rational functions then it has 
	only finitely many maximal proper invariant subvarieties.
\end{theorem}

Model-theoretically, as we show in Proposition~\ref{prop:ds-wo}, \((V,\phi)\) 
having no nonconstant invariant rational functions over~\(k\) says that the 
generic type is weakly orthogonal to the fixed field.
Hence, the model-theoretic content of Theorem~\ref{thm:dme} is that a 
rational type that is both quantifier-free internal and weakly orthogonal to 
the fixed field is isolated.
It should not be surprising, at least to the model-theorist, that this 
follows rather easily from the existence of a quantifier-free definable 
binding group.

Theorem~\ref{thm:dme} can be seen as the difference-algebraic analogue of a 
theorem in differential-algebraic geometry, appearing in~\cite{diffDME}, 
about  isotrivial \(D\)-varieties.\footnote{In fact, Proposition~2.3 
of~\cite{diffDME} proves the analogous result for the more general ``compound 
isotrivial" \(D\)-varieties -- so for types \emph{analysable}, rather than 
internal, in the constants. It is likely that compound isotriviality can be 
made sense of for rational dynamics also, and the extension of 
Theorem~\ref{thm:dme} to that case would be a desirable objective of future 
work.}

Our final application has to do with rational \(\sigma\)-varieties with the property that 
some cartesian power admits a nonconstant invariant rational 
function.
We show that there is a bound on how high a cartesian power one must look at:

\begin{theorem}[Appearing as Theorem~\ref{thm:nonorthbound} and 
	Corollary~\ref{cor:nonorthbound-auto} below]\label{thm:nonorthbound-intro}
	Suppose \((V,\phi)\) is a rational \(\sigma\)-variety over $(k,\sigma)$.
	If some cartesian power of \((V,\phi)\) admits a nonconstant invariant 
	rational function then already \((V^n,\phi)\) does, where \(n=\dim V+3\).
	In the autonomous case, when \(\phi:V\dto V\) is a rational dynamical 
	system, we can take \(n=2\), regardless of \(\dim V\).
\end{theorem}

Actually, we can weaken the antecedent of this implication somewhat, to the existence of an invariant rational function on \((V\times W,\phi\times\psi)\), for some \((W,\psi)\), that is not the pullback of a rational function on \((W,\psi)\).
See the geometric formulation of Theorem~\ref{thm:nonorthbound}, below.
As we show in Corollary~\ref{nonortho-qfint-geo}, this condition on \((V,\phi)\) turns out to be equivalent to the existence of a positive-dimensional isotrivial image.

For rational dynamics (with the bound of~\(2\)) this theorem appears already 
as~\cite[Corollary~B]{bms}, using very different methods.
But the general case is new.
An additional ingredient in its proof is the truth of the Borovik-Cherlin 
Conjecture in the theory of algebraically closed fields of characteristic 
zero, established in~\cite{nmdeg} using the work of Popov~\cite{Popov2007} as 
proposed by Borovik and Cherlin in~\cite{BC2008}.
This statement bounds the degree of generic multiple transitivity of an 
algebraic group action; and we apply that bound to the binding group action of a positive-dimensional isotrivial image of \((V,\phi)\).
The differential-algebraic analogue of Theorem~\ref{thm:nonorthbound-intro}  
(which also uses the Borovik-Cherlin Conjecture in~\(\acf_0\)) comes out of 
work in~\cite{nmdeg, pperm, abred}.
This too was part of our original motivation for developing binding groups in 
the difference-algebraic context.
Let us also mention that, while the bound of \(\dim V+3\) seems quite weak 
compared to the absolute bound of~\(2\) for rational dynamics, it turns out 
to be sharp for the analogous result in differential-algebraic geometry, and 
we expect it to be sharp here too.
But that remains as yet unverified.

Theorem~\ref{thm:nonorthbound-intro} also has a model-theoretic articulation: 
a rational quantifier-free  type~\(p\) of dimension~\(d\) is 
nonorthogonal to the fixed field if and only if the Morley power 
\(p^{(d+3)}\) is not weakly orthogonal to  the fixed field.

Finally, let us mention that one of the primary motivations for this work is 
the extension, to the setting of algebraic dynamics, of the results 
in~\cite{c3} on the structure of algebraic differential equations having the 
property that any three solutions are independent.
As that work uses model-theoretic binding groups in a crucial way, it is our 
hope that the theory developed here will lead to such an extension.

\subsection{A word about characteristic}
While we have assumed characteristic zero in this introduction, and it is required for our main theorems as stated above, much of what we do in this paper does go through for arbitrary 
characteristics.
We therefore work independently of characteristic until~$\S$\ref{sect:bg}, stating explicitly when characteristic zero is required.
While positive characteristic analogues of the theorems presented here can be 
articulated, we have decided not to do so, partly because in positive 
characteristic one should not only consider isotriviality and internality 
with respect to the fixed field, but rather to the various fixed fields 
of~\(\sigma\) composed with powers of the Frobenius automorphism.
Working out a general theory of binding groups in that setting is desirable, 
but is deferred to future work.

\subsection{Plan of the paper}
We conclude the Introduction by fixing our algebraic geometric conventions.
Then, in Section~\ref{sect:qfacfa}, we discuss/review in some detail the 
various elements of the quantifier-free fragment of~\(\acfa\)  that concern 
us.
In particular, we discuss rational types, canonical bases, nonorthogonality 
to the fixed field, and internality to the fixed field, all in the 
quantifier-free setting.
In Section~\ref{sect:rd} we develop algebraic dynamics in the general 
nonautonomous context, and produce a dictionary translating between algebraic 
dynamics and model theory.
In particular, invariant rational functions and isotriviality are discussed 
at length here.
Section~\ref{sect:bg} is dedicated to the proofs of our main binding group 
theorems, namely Theorems~\ref{thm:bgrds} and~\ref{thm:bgrt}.
Finally, in Section~\ref{sect:applications}, we state and prove the 
applications we have discussed above.

\subsection{Acknowledgements}
We would like to thank Jason Bell for numerous conversations, including one in which an error to our original proof of Corollary~\ref{cor:translational} was found and corrected.
We are also grateful for the hospitality of the Fields Institute in Toronto where some of this work was done.

\subsection{Algebraic geometric conventions}
\label{sec:ag}
Here we make explicit some more or less standard notational conventions.

We will tend to work in a sufficiently saturated algebraically closed 
field~\(\U\) that serves as a universal domain for algebraic geometry, in the 
sense of Weil.
In particular, all tuples and fields are assumed to live in~\(\U\), and all 
varieties are identified with their \(\U\)-points.
In particular, algebraic varieties are automatically reduced.

Varieties are {\bf not} assumed to be irreducible.
Subvarieties are Zariski closed.

We drop the assumption of characteristic zero, asserting it explicitly when 
needed.
If the characteristic is~\(p>0\) then we denote by~\(\fr\) the Frobenius 
automorphism of~\(\U\) given by \(x\mapsto x^p\).
In characteristic~\(0\), we take \(\fr\) to be the identity.
We denote by~\(k^{\perf}\) the perfect closure of a field~\(k\), and by 
\(k^{\alg}\) the algebraic closure.

A subfield~\(k\) is a \emph{field of definition} for an affine variety~\(V\) 
(equivalently~\(V\) is \emph{over}~\(k\))
if the ideal \(I(V)\subseteq \U[x]\) of polynomials vanishing on~\(V\) has a 
set of generators with coefficients in~\(k\).
In characteristic~\(0\) this coincides with being \(\lring\)-definable 
over~\(k\), but in positive characteristic it is a stronger notion: being 
\(\lring\)-definable over~\(k\) only ensures that \(k^{\perf}\) is a field of 
definition.

Given a variety~\(V\) over~\(k\) we will tend not to distinguish notationally 
between~\(V\) and its base change to a field extension \(K\supseteq k\), 
except when confusion could arise, in which case we use \(V_K\) for the base 
extension.

We will say that a property holds of 
\emph{general} \(a_1,\dots,a_n\in V\) to mean that it holds on a Zariski dense open 
subset of $V^n$ over~\(k\).

Given an \(n\)-tuple~\(a\), and a perfect field~\(k\),
we denote by \(\loc(a/k)\) the Zariski locus of~\(a\) over~\(k\), the 
smallest (closed) subvariety of~\(\AA^n\) over~\(k\) that contains~\(a\) as a 
\(\U\)-point.

Almost always, when working with rational maps, our varieties will be irreducible.
But as this is not always the case, let us recall how one deals with rational maps in general.
A {\em rational map} between varieties, denoted by \(\phi:V\dto W\), is an equivalence class of pairs $(U,f)$ where $U$ is a Zariski dense open subset of $V$ and $f:U\to W$ is a regular morphism, and two such pairs are deemed equivalent if the morphisms agree on the intersection of their domains.
For every rational map there is a canonical representation, $(U,f)$, such that $U$ contains the domain of any other representation -- we often denote such~$U$ by 
\(\dom{\phi}\) and identify $\phi$ with~$f$.

Suppose now that $V$ and $W$ are irreducible varieties.
A rational map $\phi:V\dto W$ is {\em dominant} if its image, or rather the image of some (equivalently any) representative of~$\phi$, is Zariski dense in~$W$.
If $\phi:V\dto W$ is dominant then it makes sense to compose~$\phi$ with a rational map on~$W$.
We say that~$\phi$ is {\em generically injective} if it has some representative $(U,f)$ such that $f$ is injective.
We say that $\phi$ is {\em birational} if it has a rational inverse; that is, if there exists a dominant rational map $\psi:W\dto V$ such that $\phi\circ\psi$ and $\psi\circ\phi$ are both the identity (as rational self-maps on~$W$ and~$V$, respectively).
Birational maps are generically injective and dominant, and the converse holds in characteristic zero.
By the {\em graph} of~$\phi$, usually denoted by $\Gamma(\phi)$, we mean the Zariski closure of $\{(a, f(a)):a\in U\}$ in $V\times W$, where $(U,f)$ is some (equivalently any) representative of~$\phi$.
Note that $\Gamma(\phi)$ is an irreducible (closed) subvariety that projects dominantly onto both~\(V\) and~\(W\), and such that the projection onto~\(V\) is birational.

We will often consider algebraic families of varieties.  These will usually 
be presented as follows: we have irreducible varieties~\(V\) and~\(Z\) over a 
field~\(k\), as well as an irreducible subvariety~\(X\subseteq V\times Z\) 
over~\(k\).
For \(e\in Z\) we denote the set-theoretic fibre by
\[
X_e:=\{v\in V:(v,e)\in X\}.
\]
That is, \(X_e\) denotes the underlying reduced variety, over \(k(e)\), of 
the subscheme of \(V_{k(e)}\) given by the scheme-theoretic fibre.
In particular, \(X_e=X_{e'}\) if and only if \(X_e(L)=X_{e'}(L)\) for some 
(equivalently any) algebraically closed field extending \(k(e)\).
These fibres form a family of subvarieties of~\(V\) parameterised by~\(Z\).  
We will tend to use \(\pi_1:X\to V\) and \(\pi_2:X\to Z\) to denote the 
co-ordinate projections.  Given another such family, say \(Y\subseteq W\times 
Z\), and a rational map
\[
\xymatrix{
X\ar[dr]\ar@{-->}[rr]^{ g}&& Y\ar[dl]\\ 
&Z
}
\]
for any \(e\in\pi_2(\dom g)\), we denote by
\(g_e:X_e\dto Y_e\)
the \(k(e)\)-rational map given by \(v\mapsto \pi_1(g(v,e))\).

It is well known that any algebraic family of subvarieties admits a rational 
quotient family where the parameters are canonical -- this is essentially a 
Hilbert scheme argument with projective varieties, but we present it here as 
a consequence of elimination of imaginaries and quantifiers.
We restrict attention to characteristic zero as we will only use this fact in 
that case, and the statement in positive characteristic is slightly more 
involved (requiring precomposition with a purely inseparable map).

\begin{fact}\label{prp:ratquotient}
	Let \(k\) be a field of characteristic zero, \(V,Z\) irreducible varieties 
	over \(k\), and \(  X\subseteq{}V\times{}Z\) an irreducible subvariety 
	projecting dominantly to \(Z\). Then there is a variety \( Z_0\) over~\(k\) 
	and a dominant rational map \(\mu=\mu_X:Z\dto Z_0\) such that
	\begin{itemize}
		\item[(a)] For general \(a,b\in{}Z\), if \(X_a=X_b\) then 
			\(\mu(a)=\mu(b)\)
		\item[(b)] Universality: For any dominant rational
			map \(t:Z\dto W\) such that for general \(a,b\in{}Z\),  \(  X_a=  X_b\) 
			implies \(t(a)=t(b)\), there is a unique
			dominant rational map \(\bar{t}: Z_0\dto W\) with \(t=\bar{t}\circ\mu\)
	\end{itemize}
	Furthermore, for general \(a,b\in{}Z\), if \(\mu(a)=\mu(b)\), then 
	\(X_a=X_b\).
\end{fact}

\begin{proof}
	By elimination of imaginaries in \(\acf\), the analogous statement holds in 
	the definable category: there is an \(\lring\)-definable map~\(\tilde\mu\) 
	on~\(Z\) such that, for all \(a,b\in{}Z\), \(X_a=X_b\) if and only if 
	\(\tilde\mu(a)=\tilde\mu(b)\).
	By quantifier elimination, and the fact that we are in characteristic zero,
	there is a nonempty Zariski open subset of~\(Z\) on which~\(\tilde\mu\) 
	agrees with a dominant rational map, \(\mu:Z\dto Z_0\).
	In particular, \(\mu\) satisfies~(a) and the ``furthermore'' clause.
	
	We claim it also satisfies~(b).
	If \(t:Z\dto W\) is dominant rational, and satisfies the condition on the 
	fibres, then there is an \(\lring\)-\emph{definable} 
	\(\tilde{t}:Z_0\ra{}W\) such that~\(t\) agrees with 
	\(\tilde{t}\circ\tilde\mu\) on a nonempty Zariski open set.
	Once again, we have \(\tilde{t}\) agrees on a nonempty Zariski open set 
	with a (necessarily dominant) rational map \(\bar{t}:Z_0\ra{}W\).
	Hence \(t=\bar{t}\circ\mu\) as rational maps on~\(Z\).
	Uniqueness is by dominance of~\(\mu\).
\end{proof}

\bigskip
\section{The quantifier-free model theory of~\(\acfa\)}\label{sect:qfacfa}

\noindent
Everything we do in this section is known to the experts, and much of it can 
be found in, or easily deduced from, the literature on the model theory of 
difference fields, in particular~\cite{acfa} and\cite{Moshe}.
Our purpose here is to give a self-contained and complete account.

Let \(\lring=\{0,1,+,-,\times\}\) be the language of rings and 
\(\lsigma=\{0,1,+,-,\times,\sigma\}\) the language of difference rings.
Fix a sufficiently saturated model \((\U,\sigma)\models\acfa\).
When working in this universal domain we follow the usual conventions that 
sets of parameters are small in cardinality compared to the degree of 
saturation -- unless explicitly stated otherwise.

We adopt the common (if unfortunate) model-theoretic convention of writing $AB$ instead of $A\cup B$ when referring to sets of parameters, except when confusion could arise.
In the same vein, if $a=(a_1,\dots,a_n)$ then we will write $Aa$ for $A\cup\{a_1,\dots,a_n\}$.

We set \(\fix(\sigma)=\{a\in\U:\sigma(a)=a\}\) to be the fixed field of 
\((\U,\sigma)\).

We do not assume, unless explicitly stated otherwise, that our difference fields are {\em inversive} -- that is they are fields~\(k\) equipped with an endomorphism~\(\sigma\), that need not be surjective.
We view difference fields as 
\(\lsigma\)-structures.
We use \(\langle A\rangle\) to denote the difference field generated by the 
set~\(A\).
If \(k\) is a  difference subfield of~\(\U\), and \(a\) is a tuple, then by 
\(k\langle a\rangle\) we mean the difference subfield of~\(\U\) generated by 
\(a\) over \(k\), namely \(k(a,\sigma(a),\sigma^2(a),\dots)\).
For natural \(m\), we denote by \(\prolong[m]{a}\) the tuple 
\((a,\sigma(a),\dots,\sigma^m(a))\).

We use nonforking independence in~\(\acfa\) freely: \(A\ind_CB\) means that 
\(\langle{}A\cup{}C\rangle^{\alg}\) is algebraically disjoint from 
\(\langle{}B\cup{}C\rangle^{\alg}\) over \(\langle{}C\rangle^{\alg}\).
In particular, dependence is always witnessed by quantifier-free formulas.

We are concerned with the quantifier-free fragment of \((\U,\sigma)\).
Note that, for substructures~$M,N$, we have that $\qftp(M/A)=\qftp(N/A)$ if and only if $M$ and $N$ are isomorphic over~$A$, and that a choice of isomorphism corresponds to a choice of enumeration.

Given a parameter set~\(A\), we will be primarily interested 
in \(S_{\qf}(A)\), the set of complete quantifier-free types over~\(A\).

\begin{definition}\label{def:qft}
	Suppose \(p\in S_{\qf}(A)\).
	We say that~\(p\) is \emph{stationary} if for any extension of parameters 
	\(B\supseteq{}A\) there is a unique extension of \(p\) to a complete 
	quantifier-free type over~\(B\) whose realisations are independent of~\(B\) 
	over~\(A\).
	This extension is the \emph{nonforking extension of~\(p\) to~\(B\)}.
\end{definition}

Complete quantifier-free types over algebraically closed difference fields 
are always stationary. In fact, let us record for future use the following 
well known strengthening of quantifier-free stationarity over algebraically 
closed sets:

\begin{lemma}\label{stationarity}
	Suppose \(k\) is an algebraically closed difference field, with four 
	difference field extensions \(K_1,K_2,L_1,L_2\).
	Assume that
	\begin{itemize}
		\item[(i)]
			\(\qftp(K_1/k)=\qftp(K_2/k)\),
		\item[(ii)]
			\(\qftp(L_1/k)=\qftp(L_2/k)\), and
		\item[(iii)]
			\(L_i\ind_kK_i\), for \(i=1,2\).
	\end{itemize}
	Then \(\qftp(K_1L_1/k)=\qftp(K_2L_2/k)\).
\end{lemma}

\begin{proof}
	Because $k$ is algebraically closed, $L_i\ind_kK$ implies that 
	\(L_i\) is linearly disjoint from \(K_i\) over \(k\), see~\cite[Lemma~2.6.7]{fried-jarden}.
	So, we have a 
	canonical identification of the field \(K_iL_i\) with the fraction field of 
	\(K_i\otimes_kL_i\).
	Moreover, this is an identification of difference fields where \(\sigma\) 
	acts on \(K_i\otimes_kL_i\) by \(\sigma(a\otimes 
	b)=\sigma(a)\otimes\sigma(b)\).
	Now \(\qftp(K_1/k)=\qftp(K_2/k)\) is witnessed by an isomorphism 
	\(\alpha:K_1\to K_2\) of difference fields over \(k\), and similarly we 
	have \(\beta:L_1\to L_2\).
	We obtain an isomorphism \(\alpha\otimes\beta:K_1\otimes_kL_1\to 
	K_2\otimes_kL_2\) of difference rings over \(k\), which will extend to the 
	fraction fields.
	That isomorphism witnesses \(\qftp(K_1L_1/k)=\qftp(K_2L_2/k)\).
\end{proof}

In particular, if \(p\in S_{\qf}(A)\) is stationary then, for all 
	\(n\geq{}1\), all \(n\)-tuples of independent realisations of~\(p\) will 
	have the same complete quantifier-free type over~\(A\), which we denote by 
	\(p^{(n)}\), and call the \emph{\(n\)-th Morley power of~\(p\)}.

Fix a difference subfield~\(k\subset\U\).

\medskip
\subsection{Rational types}\label{sec:rtypes}
Among the complete quantifier-free types over~\(k\) we will be primarily 
interested in what we will call \emph{rational types}, namely those \(p(x)\in 
S_{\qf}(k)\) that imply the formula \(\sigma(x)=f(x)\) for some rational 
function \(f\in k(x)\).
In that case, \(p\) is determined by this formula along with the 
\(\lring\)-formulas in~\(p\).

Note that if~\(p\) is rational and \(a\models p\) then 
\(k\langle{}a\rangle=k(a)\) is a finitely generated field extension of~\(k\).
Conversely, if \(k\langle a\rangle\) is finitely generated over~\(k\) as a 
field then \(\qftp(\nabla_m(a)/k)\) is rational for some \(m\geq 0\).
So, to study rational types is to study difference field extensions that are 
finitely generated as field extensions.

Here are two straightforward consequences of working with rational types that we will make use of.

\begin{lemma}
\label{lem:invalg}
If~\(k\) is inversive and \(\qftp(a/k)\) is rational then $\acl(ka)=k(a)^{\alg}$.
\end{lemma}

\begin{proof}
First observe that $k(a)^{\alg}$ is inversive.
This is seen as follows:
First, by inversiveness of~$k$, $\sigma$ extends to an isomorphism between $k(a)$ and $k(\sigma(a))$, and hence these two fields have the same transcendence degree over~$k$.
But rationality implies that $\sigma(a)\in k(a)$.
Hence $a\in k(\sigma(a))^{\alg}$.
Applying $\sigma^{-1}$ we have that $\sigma^{-1}(a)\in k(a)^{\alg}$.
Iterating, we see that $k(a)^{\alg}$ is inversive.

Since $\acl(ka)$ is the field-theoretic algebraic closure of the inversive difference field generated by $k\cup\{a\}$, it follows that $\acl(ka)=k(a)^{\alg}$.
\end{proof}

\begin{lemma}\label{lem:fgperf}
	Suppose~\(k\) is perfect and \(\qftp(a/k)\) is rational.
	If \(e\in {k(a)}^{\perf}\) then \(\qftp(\nabla_m(e)/k)\) is rational for 
	some \(m\geq0\).
\end{lemma}

\begin{proof}
	Let \(\ell>0\) be such that \(\fr^\ell(e)\in k(a)=k\langle a\rangle\).
	As a subextension of a finitely generated field extension is itself finitely generated, see~\cite[Theorem~24.9]{isaacs}, 
	it follows that \(k\langle \fr^\ell(e)\rangle\) is a finitely generated 
	field extension of~\(k\).
	Hence, for some \(m\geq 0\), if
	\[b:=\nabla_m\fr^\ell(e)=\fr^\ell\nabla_m(e)\]
	then \(\qftp(b/k)\) is rational.
	Let \(f\in k(x)\) be such that \(\sigma(b)=f(b)\).
	Then, \begin{eqnarray*}
		\sigma(\nabla_m(e))
		&=&
		\sigma(\fr^{-\ell}(b))\\
		&=&
		\fr^{-\ell}(\sigma(b))\\
		&=&
		\fr^{-\ell}(f(b))\\
		&=&
		f^{\fr^{-\ell}}(\fr^{-\ell}(b))\\
		&=&
		f^{\fr^{-\ell}}(\nabla_m(e)).
	\end{eqnarray*}
	Here \(f^{\fr^{-\ell}}\) denotes the transform of~\(f\in k(x)\) obtained by 
	applying \(\fr^{-\ell}\) to~\(k\).
	As~\(k\) is perfect, \(f^{\fr^{-\ell}}\) is again a rational function 
	over~\(k\), witnessing the rationality of \(\qftp(\nabla_m(e)/k)\).
\end{proof}

We define the \emph{dimension} of a rational type~\(p\in S_{\qf}(k)\) to be 
the transcendence degree of \(k(a)\) over \(k\) for any \(a\models p\). 

By a \emph{rational map} \(\gamma:p\to q\), between rational types \(p,q\in 
S_{\qf}(k)\), we mean that \(\gamma\) is a rational map over~\(k\) and that 
for every (equivalently some) \(a\models p\), \(\gamma(a)\models q\).
This is equivalent to asking that there are \(a\models p\) and \(b\models q\) 
with \(b\in k(a)\).
We say that~\(p\) and~\(q\) are \emph{birationally equivalent} if there exist 
rational maps \(\gamma:p\to q\) and \(\delta:q\to p\) such that \(\delta\gamma(a)=a\) and \(\gamma\delta(b)=b\) for all (equivalently some) \(a\models p\) and \( b\models q\).

The following fact about the interaction between rational types and the fixed 
field is essentially a special case of~\cite[Prop.~26]{Moshe}, but we recall 
the proof for convenience.

\begin{proposition}\label{prop:isolateC}
	Suppose~\(k\) is inversive, \(p\in S_{\qf}(k)\)  is rational, and 
	\(\C=\fix(\sigma)\).
	For any $a\models p$, the quanitifer-free type of~$a$ over $k\cup(k(a)\cap\C)$ isolates a complete quantifier-free type over $k\cup \C$.
\end{proposition}

\begin{proof}
	In what follows we use \(\tp^-\) to denote the \(\lring\)-type.

	Fix \(a_1,a_2\models p\) such that \(D:=k(a_1)\cap\C=k(a_2)\cap\C\), and 
	\(\qftp(a_1/kD)=\qftp(a_2/kD)\).
	We wish to show that \(\qftp(a_1/k\C)=\qftp(a_2/k\C)\).

	Note that the \(\sigma\)-transforms of \(a_i\) are all given by 
	\(k\)-rational functions as~\(p\) is a rational type.
	It follows that \((\tp^-(a_i/k\C)\cup{}p)\vdash\qftp(a_i/k\C)\).
	So it suffices to prove that \(\tp^-(a_1/k\C)=\tp^-(a_2/k\C)\).
	
	Let \(F_i\) be the canonical base for \(\tp^-(a_i/k\C)\) in the sense 
	of~\(\acf\).
	That is, \(F_i\) is the minimal field of definition
	of the Zariski locus of~\(a_i\) over~\(k\C\).
	Then $kF_i=k(c_i)$ for some finite tuple $c_i$ from $\C$.
	We claim that \(c_i\in\dcl(ka_i)\).
	Indeed, \(\C\) is invariant under any automorphism of \((\U,\sigma)\), and 
	hence if \(\alpha\in\aut_k(\U,\sigma)\) is such that \(\alpha(a_i)=a_i\) 
	then~\(\alpha\) preserves the Zariski locus of \(a_i\) over~\(k\C\), and 
	hence is the identity on \(F_i\), and hence on $c_i$.
	Applying Lemma~\ref{lem:invalg}, we thus have	\[
	c_i\subset\dcl(ka_i)\cap\C\subseteq{k(a_i)}^{\alg}\cap\C=:E_i.
	\]
	We claim that $E_1=E_2$.
	Indeed, given \(e\in{}E_i\), let \(P(x)\) be its minimal polynomial over 
	\(k(a_i)\).  Applying~\(\sigma\), we see that \(e\), being a tuple from~$\C$, is also a root of 
	\(P^\sigma\).
	Hence \(P=P^\sigma\) (noting that \(k(a_i)\) is a 
	\(\sigma\)-field by rationality of~\(p\)). So \(P\) is over \(\C\), and 
	hence over \(D\).
	This shows that
	\(E_i=D^{\alg}\cap\C=:E\) is independent of \(i\).
	
	It suffices to show that $\tp^-(a_1/kE)=\tp^-(a_2/kE)$.
	This is because $kE$ contains $k(c_i)=F_i$, and hence $\tp^-(a_i/k\C)$ is the unique nonforking extension of $\tp^-(a_i/kE)$.
	
	First observe that whenever~\(e\) is a finite tuple 
	from~\(E\), \(\tp^-(e/D)\vdash\tp^-(e/k(a_i))\) for \(i=1,2\).
	Indeed, let~\(\Sigma\) be the (finite) set of realisations of 
	\(\tp^-(e/k(a_1))\).
	Since \(\sigma(e)=e\) and \(\sigma(k(a_1))\subseteq k(a_1)\), we have that 
	\({\tp^-(e/k(a_1))}^\sigma\subseteq\tp^-(e/k(a_1))\), and hence 
	\(\Sigma\subseteq \sigma(\Sigma)\), which by finiteness forces 
	\(\Sigma=\sigma(\Sigma)\).
	This means that~\(\Sigma\) is \(\lring\)-definable over \(k(a_1)\cap\C=D\).
	So \(\tp^-(e/D)\vdash\tp^-(e/k(a_1))\), and similarly 
	\(\tp^-(e/D)\vdash\tp^-(e/k(a_2))\).
	
	Finally, let us show that \(\tp^-(a_1/kE)=\tp^-(a_2/kE)\).
	Given a finite tuple~\(e\) from~\(E\), we show that there is a 
	field-automorphism fixing~$k$ and taking \((a_1,e)\) to \((a_2,e)\).
	This will suffice.
	Since, by assumption, \(\qftp(a_1/kD)=\qftp(a_2/kD)\). there is a field-automorphism~$\tau$ of~\(\U\) fixing~\(kD\) pointwise and 
	such that \(\tau(a_1)=a_2\).
	Then, by the previous paragraph, \(\tp^-(\tau(e)/k(a_2))=\tp^-(e/k(a_2))\), 
	witnessed, say, by a field-automorphism~\(\iota\).
	Hence, \(\iota\tau\) fixes~$k$ and takes \((a_1,e)\) to 
	\((a_2,e)\), as desired.
\end{proof}

\medskip
\subsection{Canonical bases}
Given a quantifier-free type \(p=\qftp(a/k)\) over a perfect
difference field~\(k\),
the \emph{canonical base of~\(p\)} is the difference subfield of~\(k\) 
generated by the minimal fields of definition of the Zariski loci 
\(\loc(\nabla_n(a)/k)\), as \(n\geq 0\) varies.\footnote{This disagrees 
mildly with the terminology of Chatzidakis and Hrushovski 
in~\cite[\(\S\)2.13]{acfa}; they take as the canonical base the perfect 
closure of what we are calling the canonical base.}
Note that this does not depend on the realisation of~\(p\) chosen.
We will denote the canonical base by \(\cb(a/k)\) or by \(\cb(p)\).
When \(k\) is not necessarily perfect, we will still write \(\cb(a/k)\) to 
mean \(\cb(a/k^{\perf})\)

\begin{lemma}\label{lem:cb}
	Suppose~\(k\) is perfect and \(p=\qftp(a/k)\).
	\begin{itemize}
		\item[(a)]
			For any difference subfield \(L\subseteq k\), \(a\ind_Lk\) if and only 
			if \(\cb(a/k)\subseteq L^{\alg}\cap k\).
		\item[(b)]
			If~\(p\) is rational then \(\cb(p)\) is the difference field generated 
			by the minimal field of definition of \(\loc(\nabla(a)/k)\).
			That is, for rational types, one need only consider \(n=1\) in the 
			definition of canonical base.
		\item[(c)]
		Suppose~\(p\) is rational and~\(q=\qftp(a/K)\), where \(K\supseteq k\) is a perfect difference field extension.
			Then there is an~\(\ell\geq 0\) such that 
			\(\cb(q)\) is contained in the perfect closure of the field generated over~\(k\)
			by any~\(\ell\) independent realisations of~\(q\).
	\end{itemize}
\end{lemma}

\begin{proof}
	For part~(a) we note that
	\begin{eqnarray*}
		a\ind_Lk
		&\iff&
		\trdeg(\nabla_n(a)/L)=\trdeg(\nabla_n(a)/k),\text{ for all }n\\
		&\iff&
		\loc(\nabla_n(a)/k)\text{ is over }L^{\alg}\cap k,\text{ for all }n\\
		&\iff&
		\cb(a/k)\subseteq L^{\alg}\cap k.
	\end{eqnarray*}

	For part~(b), by rationality, we have \(\sigma(a)=f(a)\) for some rational 
	function~\(f\) over~\(k\).
	Let \(V=\loc(a/k)\).  Note that \(\Gamma:=\loc(\nabla_1(a)/k)\subseteq 
	V\times V^\sigma\) is the graph of~\(f\) viewed as a rational map on~\(V\).
	Let~\(F\) be the minimal field of definition of~\(\Gamma\).
	One shows, inductively, that \(\loc(\nabla_n(a)/k)\) is over~\(\langle 
	F\rangle\).
	Consider \(n=2\).
	Then
	\[\loc(\nabla_2(a)/k)=\loc(a,f(a),f^\sigma(f(a))/k)=\Gamma(f,f^\sigma\circ f)\]
	where $\Gamma(f,f^\sigma\circ f)$ is the graph of the rational map
	$(f,f^\sigma\circ f):V\dto V^\sigma\times V^\sigma$.
	But note that, as $\Gamma=\Gamma(f)$ and $\Gamma^\sigma=\Gamma(f^\sigma)$, we have that 
	$\Gamma(f,f^\sigma\circ f)=\Gamma\times_{V^\sigma}\Gamma^\sigma$.
	Since \(V,V^\sigma\), and \(\Gamma\) are over~\(F\), and \(\Gamma^\sigma\) 
	is over \(\sigma(F)\), we get that \(\loc(\nabla_2(a)/k)\) is 
	over~\(\langle F\rangle\).
	This argument can be iterated.

	Finally, for part~(c), let~\(F\) be the minimal field of definition of 
	\(\Gamma:=\loc(\nabla(a)/K)\).
	By a general property of canonical bases in stable theories (see~\cite[Lemma~1.2.28]{gst}), here applied to \(\acf\), there exists~\(\ell\geq 0\) such that~\(F\) is contained in the 
	perfect closure of any~\(\ell\) independent Zariski generic points of 
	\(\Gamma\).
	Let \(a_1,\dots,a_\ell\) be independent realisation of~\(q\).
	So~\(F\) is contained in the perfect closure of the field generated by 
	\(\nabla(a_1),\dots,\nabla(a_\ell)\), and hence
	\begin{eqnarray*}
		\cb(q)
		&\subseteq&\langle F\rangle\ \ \text{ by part~(b) applied to~\(q\), which is rational}\\
		&\subseteq& \langle \nabla(a_1),\dots,\nabla(a_\ell)\rangle^{\perf}\\
		&\subseteq&k(a_1,\dots,a_\ell)^{\perf}
	\end{eqnarray*}
	where the last containment uses that~\(q|_k=p\) is rational.
\end{proof}

By part~(b) of the above lemma, we have that, in the rational case, 
\(\cb(p)\) is finitely generated as a difference field.
In this case we may abuse notation by writing that \(e=\cb(p)\) to mean that 
\(\langle e\rangle\), the difference field generated by~\(e\),  is the 
canonical base of~\(p\).

\medskip
\subsection{Nonorthogonality to the fixed field}
Recall that a complete type \(\tp(a/k)\) is \emph{weakly orthogonal} to a 
\(k\)-definable set~\(\C\) if \(a\ind_kc\) for any finite tuple~\(c\) 
from~\(\C\), and it is \emph{orthogonal}\footnote{Maybe \emph{foreign} and 
\emph{weakly foreign}, as in~\cite{gst}, are better terms than orthogonal and 
weakly orthogonal, as the latter are often, and were originally, used for a 
related but symmetric notion between complete types.}
to~\(\C\) if every nonforking extension is weakly orthogonal to~\(\C\).

\begin{proposition}\label{prop:wo}
	Suppose~\(k\) is inversive, \(\C=\fix(\sigma)\), and 
	\(\qftp(a/k)\) is rational.
	Then \(\tp(a/k)\) is weakly orthogonal to~\(\C\) if and only if 
	\(k(a)\cap\C\subseteq k^{\alg}\).

	In particular, weak orthogonality to~\(\C\) depends only on \(\qftp(a/k)\).
\end{proposition}

\begin{proof}
	This is a corollary of Proposition~\ref{prop:isolateC}.

	Note that \(\tp(a/k)\) is weakly orthogonal to~\(\C\) if and only if 
	\(\tp(a/k^{\alg})\) is.
	Also, as the fixed field of an algebraic difference-field extension is 
	algebraic over the fixed field of the base, we also have that 
	\(k(a)\cap\C\subseteq k^{\alg}\) if and only if 
	\(k^{\alg}(a)\cap\C\subseteq k^{\alg}\).
	So, replacing~\(k\) by \(k^{\alg}\), we may assume that~\(k\) is 
	algebraically closed.

	The left-to-right implication is clear.
	For the converse, suppose \(k(a)\cap\C\subseteq k\) and let \(c\) be a 
	finite tuple from~\(\C\).
	Applying Proposition~\ref{prop:isolateC} to \(p:=\qftp(a/k)\), we deduce 
	that \(\qftp(a/k)\vdash\qftp(a/kc)\).
	It follows by the existence of nonforking extensions, and the 
	quantifier-free nature of nonforking, that \(a\ind_kc\).
\end{proof}

It therefore makes sense to talk about orthogonality to the fixed field for 
rational types \(p\in S_{\qf}(k)\).
Namely, \(p\) is \emph{weakly orthogonal to \(\fix(\sigma)\)} if \(a\ind_kc\) 
for some (equivalently any) \(a\models p\) and any finite tuple~\(c\) from 
\(\fix(\sigma)\); and~\(p\) is \emph{orthogonal to \(\fix(\sigma)\)} if every 
nonforking extension is weakly orthogonal to \(\fix(\sigma)\).
It turns out that to verify nonorthogonality one need not consider all 
nonforking extensions:

\begin{proposition}\label{o-wo}
	Suppose~$k$ is inversive and algebraically closed, and \(p\in S_{\qf}(k)\) is rational.
	Then~\(p\) is nonorthogonal to \(\fix(\sigma)\) if and only if 
	\(p^{(\ell)}\) is not weakly orthogonal to \(\fix(\sigma)\), for some 
	\(\ell\geq 1\).
\end{proposition}

\begin{proof}
	This is a standard argument using canonical bases and forking calculus.

	The right-to-left direction is clear.
	For the converse, suppose~\(p\) is nonorthogonal to \(\fix(\sigma)\), and 
	let this be witnessed by a difference field extension \(K\supseteq k\), 
	\(a\models p\) with \(a\ind_kK\), and~\(c\) from~\(\fix(\sigma)\) with 
	\(a\nind_Kc\).
	Extending~\(K\), if necessary, we may assume that~\(K\) is perfect.
	Let~\(e=\cb(ac/K)\).
	Since \(\qftp(ac/K)\) is rational, Lemma~\ref{lem:cb}(c) gives us that
	there are independent realisations \(a_1c_1,a_2c_2,\dots,a_\ell c_\ell\) of 
	\(\qftp(ac/K)\) such that \(e\in {k(a_1c_1,\dots,a_\ell c_\ell)}^{\perf}\).
	(Here~\(\ell\) could be~\(0\).)
	Moreover, we can choose the \(a_ic_i\) such that \(ac\ind_Ka_1c_1\dots 
	a_\ell c_\ell\).

	We first claim that
	\(a\nind_{ke} c\).
	Indeed, from \(a\nind_Kc\) we get that
	\begin{equation}\label{nind}
		a\nind_{ke}Kc,
	\end{equation}
	and from \(ac\ind_{ke}K\) we get that
	\begin{equation}\label{ind}
		a\ind_{kec}K.
	\end{equation}
	From~(\ref{nind}) and~(\ref{ind}) we get the desired \(a\nind_{ke} c\).

	Next, we claim that
	\(a\nind_{ka_1c_1\dots a_\ell c_\ell}c\).
	Indeed, from \(ac\ind_{ke}K\) and \[ac\ind_Ka_1c_1\dots a_\ell c_\ell\] we 
	get that
	\(ac\ind_{ke}Ka_1c_1\dots a_\ell c_\ell\).
	In particular,
	\(a\ind_{ke}a_1c_1\dots a_\ell c_\ell\).
	So, if it were the case that \(a\ind_{ka_1c_1\dots a_\ell c_\ell}c\) then 
	we would have
	\(a\ind_{ke}ca_1c_1\dots a_\ell c_\ell\)
	which contradicts
	\(a\nind_{ke}c\).
	(Here we are using that
	\(e\in {k(a_1c_1,\dots, a_\ell c_\ell)}^{\alg}\)
	in order to apply the transitivity of nonforking.)

	Finally, from \(a\nind_{ka_1c_1\dots a_\ell c_\ell}c\) it follows that
	\((a,a_1,\dots,a_\ell)\nind_k(c,c_1,\dots,c_\ell)\).
	This suffices as \((a,a_1,\dots,a_\ell)\models p^{(\ell+1)}\) and 
	\((c,c_1,\dots,c_\ell)\) is a tuple from the fixed field.
\end{proof}

One of our main applications of binding groups is the existence of a bound 
on~\(\ell\) in the statement of the above proposition -- this is 
Theorem~\ref{thm:nonorthbound} below.

We will make use of the following immediate corollary, which says that 
nonorthogonality to the fixed field is always witnessed by parameters that 
themselves realise rational types:

\begin{corollary}\label{cor:o-wo}
	Suppose~\(k\) is inversive and algebraically closed, and \(p=\qftp(a/k)\) is a rational 
	type that is nonorthogonal to the fixed field.
	Then there is a tuple~\(b\) such that \(\qftp(b/k)\) is rational, 
	\(a\ind_kb\), and \(a\nind_{kb}c\) for some tuple~\(c\) from 
	\(\fix(\sigma)\).
\end{corollary}

\begin{proof}
	Let~\(\ell\geq 1\) be such that \(p^{(\ell)}\) is not weakly orthogonal to 
	the fixed field, and let \((a=a_1,\dots,a_\ell)\models p^{(\ell)}\).
	Then \(b=(a_2,\dots,a_\ell)\) has the desired properties.
\end{proof}

\medskip
\subsection{Quantifier-free internality to the fixed 
field}\label{subsect:qfint}
While orthogonality to the fixed field behaves well with the quantifier-free 
fragment of \(\acfa\), at least for rational types, \emph{internality} is 
harder to pin down because we do not quite understand \(\dcl\) in \(\acfa\).
Indeed, observe that, in general, definable closure properly contains quantifier-free definable closure:
assuming that the square roots of $-1$ are not fixed by~$\sigma$, every nonzero $a\in \fix(\sigma)$ is definable from $a^2$ but not quantifier-free definable.
Moreover, quantifier-free definable closure properly contains perfect closure: we can have $a\in\fix(\sigma)$ such that $a$ is the only one of the $3$rd roots of $a^3$ that is in $\fix(\sigma)$, so that $a$ is quantifier-free definable from $a^3$ without being in its perfect closure (as long as the characteristic is not~$3$).
A detailed analysis of definable closure in the fixed field can be found in~\cite{beyarslan}.

Following~\cite[\(\S\)5]{acfa}, we therefore take a rather strong condition 
for our notion of quantifier-free internality:

\begin{definition}[Quantifier-free internality]\label{def:type-qfint}
	Suppose~\(\C\) is a quantifier-free 
	\(k\)-definable set and \(p\in S_{\qf}(k)\) is stationary.
	We say that \(p\) is \emph{qf-internal to \(\C\)} if for all \(a\models p\) 
	there is a difference field extension \(K\supseteq k\) such that 
	\(a\ind_kK\) and
	\(a\in K\langle c\rangle^{\perf}\) for some tuple \(c\) from \(\C\).
\end{definition}

The condition is strong in that we ask for~\(a\) to be in the perfect closure 
of \(K\langle c\rangle\) rather than simply to be quantifier-free definable 
from~\(c\) over~\(K\).
In fact, as the following proposition shows, this condition is even stronger 
than it looks when we restrict to rational types and the fixed field:

\begin{proposition}\label{prop:equivalences-qfint}
	Suppose~$k$ is algebraically closed and \(p\in S_{\qf}(k)\) is rational.
	Then the following are equivalent:
	\begin{itemize}
		\item[(i)]
			\(p\) is qf-internal to \(\fix(\sigma)\).
		\item[(ii)]
			For all (equivalently for some) \(a\models p\) there is a difference 
			field extension \(K\supseteq k\) such that \(a\ind_kK\) and \(a\in 
			K(c)\) for some tuple~\(c\) from \(\fix(\sigma)\).
		\item[(iii)]
			For all (equivalently for some) \(a\models p\) there is a difference 
			field extension \(K\supseteq k\) and \(c\) from~\(\fix(\sigma)\) such 
			that \(a\ind_kK\) and \(K(a)\subseteq K(c)\subseteq{K(a)}^{\perf}\).
		\item[(iv)]
			For all \(a\models p\) there exists \(K=k(a_1,\dots,a_n,d)\) where
			\begin{itemize}
				\item
					\(a_1,\dots,a_n\) are independent realisations of \(p\) over \(k\), 
					and
				\item
					\(d\) is a tuple from~\(\fix(\sigma)\),
			\end{itemize}
			such that
			\(a\ind_kK\) and \({K(a)}^{\perf}={K(c)}^{\perf}\), for some~\(c\) 
			from~\(\fix(\sigma)\).
	\end{itemize}
\end{proposition}

\begin{proof}
	Let \(\C:=\fix(\sigma)\).

	Assuming (i) we prove the ``for all \(a\models p\)" version of (ii).
	Let \(a\models p\).
	By definition, there is a difference field extension \(K\supseteq k\) such 
	that
	\(a\ind_kK\) and \(a\in K\langle c\rangle^{\perf}\) where \(c\) is a tuple 
	from \(\C\).
	We may assume that~\(K\) is perfect.
	As \(c\) is in the fixed field, we have that \(a\in{K(c)}^{\perf}\).
	It follows that for some \(\ell\geq 0\) and \(c':=\fr^{-\ell}(c)\), we have 
	\(a\in K(c')\).
	As \(\C\) is perfect, \(c'\) is also from \(\C\).

	Assuming the ``for some \(a\models p\)" version of (ii)  we prove the ``for 
	some \(a\models p\)" version of~(iii).
	Let \(a\models p\) and \(K\supseteq k\) be such that \(a\ind_kK\) and 
	\(a\in K(c)\), where~\(c\) is a tuple from \(\C\).
	Replacing \(K\) by \(K(c_0)\) where \(c_0\) is a maximal sub-tuple of \(c\) 
	such that \(a\ind_kKc_0\), we may assume that \(c\in{K(a)}^{\alg}\).
	We may also assume that~\(K\) is perfect.
	Let \(e=\cb(c/K(a))\in{K(a)}^{\perf}\).
	In terms of the pure algebraically closed field~\(\U\), this means 
	that~\(e\) is a code for the (finite) set~\(E\) of \(Ka\)-conjugates 
	of~\(c\).
	As~\(\sigma\) fixes~\(c\), it fixes~\(E\), and hence~\(e\) is a tuple 
	from~\(\C\).
	On the other hand, as \(a\in K(c)\) we have that \(a=f(c)\) for some~\(f\) 
	a rational function over~\(K\).
	So \(a=f(c')\) for \emph{every} \(c'\in E\).
	Hence, any field automorphism fixing \(K(e)\) will fix \(a\), proving that 
	\(a\in{K(e)}^{\perf}\).
	Replacing~\(e\) by some \(e'=\fr^{-\ell}(e')\) we have that \(a\in K(e')\), 
	and it is still the case that \(e'\) is from~\(\C\) and that 
	\(e'\in{K(a)}^{\perf}\).

	Next, we show that the ``for some" version of~(iii) implies the ``for all" 
	version.
	Fix \(a\models p\) and \(K\supseteq k\) satisfying~(iii).
	Let \(a'\models p\) be another realisation.
	Choose \(K'\models\tp(K/ka)\) with \(K'\ind_{ka}a'\).
	Then \(K'\ind_kaa'\).
	By Lemma~\ref{stationarity}, \(\qftp(Ka/k)=\qftp(K'a'/k)\).
	Hence, \(K(a)\) and \(K'(a')\) are difference-field isomorphic over~\(k\).  
	The fact that there is \(c\) from \(\C\) such that \(K(a)\subseteq 
	K(c)\subseteq{K(a)}^{\perf}\) implies that there must be some \(c'\) from 
	\(\C\), namely the image of \(c\) under the above isomorphism, such that 
	\(K'(a')\subseteq K'(c')\subseteq {K'(a')}^{\perf}\), as desired.

	(iii)\(\implies\)(iv).
	This is similar to the proof of Proposition~\ref{o-wo}.
	Fix \(a\models p\), and let \(K\supseteq k\) and~\(c\) 
	from~\(\fix(\sigma)\) satisfying~(iii).
	We may assume that~\(K\) is perfect.
	Let \(e=\cb(ac/K)\).
	The fact that \({K(a)}^{\perf}={K(c)}^{\perf}\) is reflected in 
	\(\loc(a,c/K)\), whose minimal field of definition is contained 
	in~\(k\langle e\rangle\).
	It follows that \({k\langle 
	e\rangle(a)}^{\perf}={k\langle{}e\rangle(c)}^{\perf}\).
	On the other hand, Lemma~\ref{lem:cb}(c) gives us that
	there are independent realisations \(a_1c_1,a_2c_2,\dots,a_\ell c_\ell\) of 
	\(\qftp(ac/K)\) such that \(e\in{k(a_1c_1,\dots,a_\ell c_\ell)}^{\perf}\).
	Moreover, we can choose the \(a_ic_i\) such that \(a\ind_Ka_1c_1\dots 
	a_\ell c_\ell\) and hence \(a\ind_ka_1c_1\dots a_\ell c_\ell\).
	Hence,
	letting \(K':=k(a_1c_1,\dots,a_\ell c_\ell)\), we get that \(a\ind_kK'\) 
	and
	\({K'(a)}^{\perf}={K'(c)}^{\perf}\).
	Finally, observe that \(K'\) is of the form called for by~(iv).

	(iv)\(\implies\)(i) is clear.
\end{proof}

The following proposition shows that qf-internality to the fixed field arises 
whenever there is nonorthogonality.

\begin{proposition}\label{nonortho-qfint}
	Suppose~$k$ is inversive and algebraically closed, and \(p\in S_{\qf}(k)\) is rational.
	The following are equivalent:
	\begin{itemize}
		\item[(i)]
			\(p\) is nonorthogonal to~\(\fix(\sigma)\).
		\item[(ii)]
			There is a rational map \(p\to q\) where \(q\in S_{\qf}(k)\) is 
			positive-dimensional rational and qf-internal to~\(\fix(\sigma)\).
	\end{itemize}
\end{proposition}

\begin{proof}
	Let \(\C=\fix(\sigma)\).

	(i)\(\implies\)(ii).
	Let \(a\models p\).
	By Corollary~\ref{cor:o-wo},
	there is~\(b\) such that \(\qftp(b/k)\) is rational, \(a\ind_kb\), and 
	\(a\nind_{kb}c\) for some tuple~\(c\) from \(\C\).
	Consider \(r:=\qftp(bc/k(a))\) and the canonical base 
	\(e:=\cb(r)\in{k(a)}^{\perf}\).
	As~\(r\) is rational and \(bc\nind_{k}a\), Lemma~\ref{lem:cb}(a) implies 
	that \(e\notin k\).
	Moreover, by part~(c) of that lemma, there are independent realisations 
	\(b_1c_1,\dots, b_nc_n\) of \(r\), such that
	\(e\in{k(b_1c_1,\dots, b_nc_n)}^{\perf}\).
 
	Since \(e\in{k(a)}^{\perf}\) and \(p\) is rational, Lemma~\ref{lem:fgperf} 
	implies that \(\qftp(\nabla_m(e)/k)\) is rational, for some \(m\geq 0\).
	There is also some \(\ell\geq 0\) such that
	\(\fr^\ell(\nabla_m(e))\in k(a)\).
	Let \(e':=\fr^\ell(\nabla_m(e))\) and set \(q:=\qftp(e'/k)\).
	Then \(q\) is rational and positive-dimensional, and we have a rational map 
	\(p\to q\).
	It remains to show that~\(q\) is qf-internal to~\(\C\).

	Let \(b':=(b_1,\dots,b_n)\).
	We claim that \(e'\ind_kb'\).
	Indeed, as \(b_1,\dots, b_n\) are independent realisations of 
	\(\qftp(b/k(a))\), and \(b\ind_ka\), it follows that \(a\ind_kb'\).
	Since \(e'\in k(a)\), we get \(e'\ind_kb'\), as desired.

	Let \(c':=(c_1,\dots,c_n)\).
	Recall that \(e\in{k(b'c')}^{\perf}\).
	Increasing \(\ell\) if necessary, we may assume that \(e'\in k(b'c')\).
	As~\(c'\) is a tuple from~\(\C\), this witnesses that~\(q\) is qf-internal 
	to~\(\C\).

	(ii)\(\implies\)(i).
	Fix \(e\models q\).
	By qf-internality to~\(\C\), there is \(K\supseteq k\), and~\(c\) 
	from~\(\C\), such that \(e\ind_kK\) and \(e\in K(c)\).
	Since \(\dim(q)>0\), we have that \(e\notin K^{\alg}\), and hence 
	\(e\nind_Kc\).
	Choose \(a\models p\) such that \(e\in k(a)\) and \(a\ind_{ke}K\).
	We get that \(a\ind_kK\) and \(a\nind_Kc\), witnessing that~\(p\) is 
	nonorthogonal to~\(\C\).
\end{proof}

\bigskip
\section{Algebraic dynamics}\label{sect:rd}

\noindent
We will see that studying rational types corresponds to a certain general 
setting for algebraic dynamics.
In this section we develop the basics of this theory.

We fix, for the remainder of this paper, a base difference field $(k,\sigma)$ that is inversive and algebraically closed.

By a \emph{rational \(\sigma\)-variety over~\(k\)} we mean an 
irreducible
variety \(V\) over~\(k\) equipped with a dominant rational map \(\phi:V\dto 
V^\sigma\) over~\(k\).
Here \(V^\sigma\) denotes the transform of~\(V\) with respect to the action 
of~\(\sigma\) on the field of definition~\(k\).
Note that the rational \(\sigma\)-variety structures on~\(V\) correspond 
precisely to the extensions of~\(\sigma\) from~\(k\) to the rational function 
field \(k(V)\), given by \(f\mapsto f^\sigma\circ\phi\).

Let us emphasise our (somewhat unfortunate) convention that while varieties need not be irreducible in general, the underlying variety of a rational \(\sigma\)-variety is assumed to be irreducible.

If~\(\sigma\) is trivial on~\(k\) then we say that \((V,\phi)\) is a
\emph{rational dynamical system};
in that case \(\phi:V\dto V\) is a rational transformation of~\(V\).
This is often the setting that algebraic dynamics is restricted to, but we 
will work generally.

The study of \(\sigma\)-varieties as a geometric category in its own right 
was initiated in~\cite{kp2007} and applied to algebraic dynamics 
in~\cite{ch1,ch2}.

We will sometimes be interested in the cartesian powers of a rational 
\(\sigma\)-variety, which we will tend to denote either by \({(V,\phi)}^n\) 
or by \((V^n,\phi)\).
We mean, of course, the rational \(\sigma\)-variety whose underlying variety 
is the cartesian power \(V^n\) equipped with the dominant rational map 
\(V^n\dto{(V^n)}^\sigma={(V^\sigma)}^n\) given co-ordinatewise by~\(\phi\), 
but which we continue to denote by \(\phi\).
Irreducibility is preserved in passing to cartesian powers (or 
products) by our assumption that~$k$ is algebraically closed.

An \emph{invariant subvariety} of \((V,\phi)\) is an irreducible Zariski closed subset~\(X\) of~\(V\) over~\(k\) such that \(X\cap\dom\phi\) is nonempty and 
\(\phi(X)\) is Zariski dense in \(X^\sigma\).
Equivalently, \((X,\phi|_X)\) is itself a  rational \(\sigma\)-variety.

By an \emph{equivariant rational map} \( g:(V,\phi)\dto (W,\psi)\) we mean 
that \( g:V\dto W\) is a rational map and that \(\psi g= g^\sigma \phi\) as 
rational maps from \(V\) to \(W^\sigma\).
The equivariant map \(g\) is said to be \emph{dominant}, \emph{birational}, 
etc., if it is such as a rational map of algebraic varieties. Note that the 
inverse of an equivariant birational map of \(\sigma\)-varieties is itself 
equivariant.

\begin{lemma}\label{lem:graph}
	Suppose \((V,\phi)\) and \((W,\psi)\) are rational 
	\(\sigma\)-varieties over \((k,\sigma)\), and \(f:(V,\phi)\dto(W,\psi)\) is 
	a dominant equivariant rational map.
	Then the graph of~\(f\) is an invariant subvariety of \((V\times 
	W,\phi\times\psi)\).
\end{lemma}

\begin{proof}
	Let \(\Gamma(f)\subseteq V\times W\) denote the graph of~\(f\).
	It is an irreducible closed subvariety that projects dominantly onto 
	both~\(V\) and~\(W\), and such that the projection onto~\(V\) is a 
	birational equivalence.

	Since~\(f:V\dto W\) is dominant it takes the nonempty Zariski open subset 
	\(\dom\phi\cap\dom f\) of~\(V\) to a Zariski dense subset of~\(W\).
	In particular, there exists a point \(v\in\dom\phi\cap\dom f\) such that 
	\(f(v)\in\dom\psi\).
	Hence, \((v, f(v))\) witnesses that \(\Gamma(f)\cap\dom{\phi\times\psi}\) 
	is nonempty.

	Next, we show that \(\phi\times\psi\) takes \(\Gamma(f)\) to 
	\({\Gamma(f)}^\sigma\).
	Fix \((v,f(v))\in\Gamma(f)\cap\dom{\phi\times\psi}\).
	Then
	\begin{eqnarray*}
		(\phi\times\psi)(v,f(v))
		&=&
		\big(\phi(v),\psi(f(v))\big)\\
		&=&
		\big(\phi(v),f^\sigma(\phi (v))\big)\ \ \ \text{ by equivariance}\\
		&\in&
		\Gamma(f^\sigma)\\
		&=&
		{\Gamma(f)}^\sigma.
	\end{eqnarray*}
	Since \(\Gamma(f)\cap\dom{\phi\times\psi}\) is Zariski dense in 
	\(\Gamma(f)\), it follows that \(\phi\times\psi\) takes all of 
	\(\Gamma(f)\) to \({\Gamma(f)}^\sigma\).

	Finally, to show Zariski-density of the image, work over any field 
	extension \(K\supseteq k\) and let \(v\in V\) be Zariski generic, so that 
	\((v,f(v))\) is Zariski generic in \(\Gamma(f)\).
	Then, by dominance of \(\phi:V\dto V^\sigma\), we have that \(\phi(v)\) is 
	Zariski generic in \(V^\sigma\).
	And so, \(\big(\phi(v),f^\sigma(\phi (v))\big)\) is Zariski generic in 
	\(\Gamma(f^\sigma)\).
	But \(\big(\phi(v),f^\sigma(\phi (v))\big)=(\phi\times\psi)(v,f(v))\) and 
	\(\Gamma(f^\sigma)={\Gamma(f)}^\sigma\), so that
	\((\phi\times\psi)(v,f(v))\) is Zariski generic in \({\Gamma(f)}^\sigma\), 
	as desired.
\end{proof}

Fix, now, a sufficiently saturated model \((\U,\sigma)\models\acfa\) 
extending \((k,\sigma)\).
Associated to a rational \(\sigma\)-variety \((V,\phi)\) is the 
quantifier-free \(\lsigma\)-definable set
\[
{(V,\phi)}^\sharp:=\{a\in\dom\phi:\sigma(a)=\phi(a)\}
\]
with parameters from~\(k\).
Existential closedness of \((\U,\sigma)\) ensures that that this set is nonempty.
In fact, it is Zariski dense in~\(V\).
Moreover, we can associate to \((V,\phi)\) a rational type \(p(x)\in 
S_{\qf}(k)\), the \emph{generic quantifier-free type of \((V,\phi)\) 
over~\(k\)}, which is determined by saying that~\(x\) is Zariski generic 
in~\(V\) over~\(k\) and \(x\in{(V,\phi)}^\sharp\).
By a \emph{generic point of \((V,\phi)\)} we mean a 
realisation of this generic type.

Every rational type arises in this way.
Indeed, given \(p\in S_{\qf}(k)\) rational, fix \(a\models p\), let 
\(V=\loc(a/k)\) be the Zariski locus of~\(a\) over~\(k\), and 
take~\(\phi:V\dto V^\sigma\) to be the rational map whose graph is 
\(\loc(a,\sigma(a)/k)\).
That this locus is the graph of a rational map is a consequence of the fact 
that~\(p\) is a rational type.
Because $(k,\sigma)$ is inversive, $\sigma(a)$ is Zariski generic in $V^\sigma$ over~$k$, and so~$\phi$ is dominant.
So~\(p\) is the generic quantifier-free type of \((V,\phi)\).

These constructions are functorial:
given rational \(\sigma\)-varieties \((V,\phi)\) and \((W,\psi)\), with 
generic quantifier-free types~\(p\) and~\(q\), respectively, dominant 
equivariant rational maps \((V,\phi)\dto(W,\psi)\) correspond (via 
restriction) to rational maps \(p\to q\) in the sense of 
Section~\ref{sec:rtypes}.

We will be considering algebraic families of \(\sigma\)-varieties, and we 
record for future use the fact that they can be made canonical (at least in 
characteristic zero):

\begin{proposition}\label{prop:canfam}
	Suppose \(\characteristic(k)=0\) and \((V,\phi)\) and \((Z,\psi)\) are rational \(\sigma\)-varieties over \((k,\sigma)\), 
	and \(\Gamma\subseteq V\times Z\) is an irreducible Zariski closed subset which is 
	invariant for \(\phi\times\psi\), and such that \(\pi_1:\Gamma\to Z\) is 
	dominant.
	There exists a rational \(\sigma\)-variety \((Z_0,\psi_0)\) and an 
	equivariant dominant rational map \(\mu:(Z,\psi)\dto(Z_0,\psi_0)\) such 
	that for general \(a,a'\in Z\),
	\(
	\mu(a)=\mu(a')\ \iff \ \Gamma_a=\Gamma_{a'}
	\).
\end{proposition}

\begin{proof}
	From \(\Gamma\subseteq V\times Z\), Fact~\ref{prp:ratquotient} provides a 
	dominant rational map \(\mu:Z\dto Z_0\) such that, for general \(a,a'\in 
	Z\),
	\(\mu(a)=\mu(a')\) if and only if \(\Gamma_a=\Gamma_{a'}\).
	It remains, therefore, to put a \(\sigma\)-variety structure on \(Z_0\) 
	such that~\(\mu\) is equivariant.

	Let \(f:Z\dto Z_0\) be \(\mu^\sigma\circ\psi\). We claim that for general 
	\(a,b\in{}Z\), if \(\Gamma_a=\Gamma_b\) then \(f(a)=f(b)\).
	Since \(\mu^\sigma=\mu_{\Gamma^\sigma}\) is a quotient map for 
	\(\Gamma^\sigma\), it suffices to show that if \(\Gamma_a=\Gamma_b\) then 
	\({\Gamma^\sigma}_{\psi(a)}={\Gamma^\sigma}_{\psi(b)}\).
	Fix~\(x\in \Gamma^\sigma_{\psi(a)}\) a Zariski generic point over 
	\(k(a,b)\).  Then \((\psi(a),x)\in\Gamma^\sigma\) is Zariski generic 
	over~\(k\), and hence, by \((\psi\times\phi)\)-invariance, is of the form 
	\((\psi(a),\phi(v))\) for some \(v\in V\) such that \((a,v)\in\Gamma\).  It 
	follows that \(v\in\Gamma_a=\Gamma_b\), and so \((b,v)\in\Gamma\) and 
	Zariski generic over~\(k\).  By  \((\psi\times\phi)\)-invariance again, 
	\(x=\phi(v)\in\Gamma^\sigma_{\psi(b)}\).  As this is the case for all 
	Zariski generic points over \(k(a,b)\), it follows that 
	\({\Gamma^\sigma}_{\psi(a)}\subseteq{\Gamma^\sigma}_{\psi(b)}\), and we 
	conclude \({\Gamma^\sigma}_{\psi(a)}={\Gamma^\sigma}_{\psi(b)}\), by 
	symmetry.

	It now follows from the universality of~\(\mu\), given by 
	Fact~\ref{prp:ratquotient}, that there is a unique dominant rational map 
	\(\psi_0:Z_0\ra{}{Z_0}^\sigma\) with 
	\(\psi_0\circ\mu=f=\mu^\sigma\circ\psi\), as required.
\end{proof}

\medskip
\subsection{Invariant rational functions}\label{subsect:irf}
Of special interest are equivariant rational maps from \((V,\phi)\) to the affine line equipped with the trivial dynamics, namely
\(\lambda:(V,\phi)\dto(\AA^1,\id)\).
These are called the \emph{invariant rational functions} on \((V,\phi)\); 
they are those rational functions, \(\lambda\in k(V)\), such that 
\(\lambda=\lambda^\sigma\circ\phi\).

\begin{lemma}\label{lem:invratchar}
	Suppose~\(a\) is a generic point of \((V,\phi)\), and \(\lambda\in k(V)\).
	Then \(\lambda\) is an invariant rational function of \((V,\phi)\) if and 
	only if  \(\lambda(a)\in\fix(\sigma)\).
\end{lemma}

\begin{proof}
	Note that
	\(\sigma(\lambda(a))=\lambda^\sigma(\sigma(a))=\lambda^\sigma(\phi(a))\),
	where the final equality is because \(a\in{(V,\phi)}^\sharp\).
	Hence, if \(\lambda=\lambda^\sigma\circ\phi\) then 
	\(\lambda(a)\in\fix(\sigma)\).
	Conversely, if \(\lambda(a)\in\fix(\sigma)\) then
	\(\lambda(a)=\lambda^\sigma(\phi(a))\).
	But, as \(a\) is Zariski generic in~\(V\), it follows that 
	\(\lambda=\lambda^\sigma\phi\) as rational functions on~\(V\).
\end{proof}

We have the following geometric characterisation of nonorthogonality to the 
fixed field in terms of invariant rational functions:

\begin{proposition}\label{prop:ds-wo}
	Suppose \((V,\phi)\) is a rational \(\sigma\)-variety over~\((k,\sigma)\) with quantifier-free generic 
	type~\(p\).
	\begin{itemize}
		\item[(a)]
			\(p\) is weakly orthogonal to \(\fix(\sigma)\) if and only if 
			\((V,\phi)\) admits no nonconstant invariant rational functions.
		\item[(b)]
			\(p\) is orthogonal to \(\fix(\sigma)\) if and only if, for every 
			rational \(\sigma\)-variety \((W,\psi)\) over~\(k\), the invariant 
			rational functions on \((V,\phi)\times(W,\psi)\) are all pullbacks of 
			invariant rational functions on \((W,\psi)\).
	\end{itemize}
\end{proposition}

\begin{proof}
	Fix $a\models p$.
	By Proposition~\ref{prop:wo}, along with the fact 
	that~\(k\) is algebraically closed, we have that~$p$ is weakly orthogonal to $\fix(\sigma)$ if and only if $k(a)\cap\fix(\sigma)\subseteq k$.
	By Lemma~\ref{lem:invratchar}, this latter condition says precisely that $(V,\phi)$ admits no nonconstant invariant rational functions, proving~(a).

	Towars part~(b), suppose that \((W,\psi)\) is another rational \(\sigma\)-variety 
	over~\(k\) and \(f\) is an invariant rational function on 
	\((V,\phi)\times(W,\psi)\).
	Let \(a\models p\) and~\(b\) generic in \((W,\psi)\), with \(a\ind_kb\).
	Then~\(a\) is generic in \((V,\phi)\) over \(K:=k(b)\), and 
	\(\lambda:=f(-,b)\) is a rational function on~\(V\) over~\(K\) with the 
	property that
	\[
	\sigma(\lambda(a))=\lambda^\sigma(\sigma(a))=
	f^\sigma(\sigma(a),\sigma(b))=\sigma(f(a,b))=f(a,b)=\lambda(a).
	\]
	That is, \(\lambda\) is an invariant rational function on \((V,\phi)\) 
	over~\(K\).
	If~\(p\) is orthogonal to \(\fix(\sigma)\) then the nonforking extension 
	of~\(p\) to~\(K\) is weakly orthogonal to \(\fix(\sigma)\), and hence, by 
	part~(a), we have that~\(\lambda\in K^{\alg}\).
	It follows from the absolute irreducibility of~\(V\) that~\(K=k(W)\) is 
	relatively algebraically closed in~\(k(V\times W)\), and hence \(\lambda\in 
	K\).
	Writing \(\lambda=g(b)\)
	we see that~\(g\) is an invariant rational function on \((W,\psi)\) and  
	\(f\) is the pullback of~\(g\).
	This proves the left-to-right implication of part~(b).

	For the converse, suppose~\(p\) is nonorthogonal to \(\fix(\sigma)\), and 
	let this be witnessed by \(B\supseteq k\) and~\(c\) from~\(\fix(\sigma)\) 
	such that \(a\ind_kB\) and \(a\nind_Bc\).
	By Corollary~~\ref{cor:o-wo}, we can choose~\(B\) of the form \(kb\) where 
	\(q:=\qftp(b/k)\) is rational.  Let \((W,\psi)\) be a rational 
	\(\sigma\)-variety over~\(k\) such that~\(q\) is the generic 
	quantifier-free type of \((W,\psi)\).
	The fact that \(a\nind_{kb}c\) tells us that \(\tp(a/K)\), where 
	\(K:=k(b)\), is not weakly orthogonal to \(\fix(\sigma)\).
	Hence, by Proposition~\ref{prop:wo}, there exists 
	\(\lambda\in{}K(a)\cap\fix(\sigma)\setminus{}K^{\alg}\).
	Writing \(\lambda=f(a,b)\) we have that \(f\) is an invariant rational 
	function on \((V,\phi)\times(W,\psi)\).
	The fact that \(\lambda\notin{}K\) tells us that~\(f\) is not the pullback 
	of a rational function on~\(W\).
\end{proof}

\medskip
\subsection{Isotriviality}
The geometric  counterpart to quantifier-free internality to the fixed field 
is isotriviality in the following natural sense:

\begin{definition}\label{def:iso}
	Suppose \((V,\phi)\) is a rational 
	\(\sigma\)-variety over~\((k,\sigma)\).
	By a \emph{trivialisation of \((V,\phi)\) over~\(k\)} we mean
	\begin{itemize}
		\item
			a rational \(\sigma\)-variety \((Z,\psi)\),
		\item
			an invariant subvariety~\(Y\) of \((\AA^\ell\times Z,\id\times\psi)\), 
			and,
		\item
			an equivariant birational map
			\[
			\xymatrix{
			(V\times Z,\phi\times\psi)\ar[dr]\ar@{-->}[rr]^{g}_{\isom}&& 
			(Y,\id\times \psi)\ar[dl]\\
			&(Z,\psi)
			}
			\]
	\end{itemize}
	all defined over~\(k\).
	We say that \((V,\phi)\) is \emph{isotrivial} if there exists a 
	trivialisation.
\end{definition}

\begin{proposition}\label{prop:ds-iso}
	Suppose \((V,\phi)\) is a rational \(\sigma\)-variety over~\((k,\sigma)\) with generic quantifier-free 
	type~\(p\).
	\begin{itemize}
		\item[(a)]
			If \((V,\phi)\) is isotrivial then~\(p\) is qf-internal to 
			\(\fix(\sigma)\).
		\item[(b)]
			Suppose \(\characteristic(k)=0\).
			If~\(p\) is qf-internal to \(\fix(\sigma)\) then \((V,\phi)\) is 
			isotrivial.
	\end{itemize}
\end{proposition}

\begin{proof}
	Suppose
	\[
	\xymatrix{
	(V\times Z,\phi\times \psi)\ar[dr]\ar@{-->}[rr]^{ g}_{\isom}&& 
	(Y,\id\times \psi)\ar[dl]\\
	&(Z,\psi)
	}
	\]
	is a trivialisation of \((V,\phi)\) over~\(k\).
	Choose~\(a\) generic in \((V,\phi)\) and~\(b\) be generic in~\((Z,\psi)\), 
	with \(a\ind_kb\).
	Hence, \((a,b)\) is generic in \((V\times Z,\phi\times \psi)\), so that 
	\(g(a,b)\) is generic in \((Y,\id\times \psi)\).
	In particular, \(g(a,b)\in{(Y,\id\times \psi)}^\sharp\), so that 
	\(g(a,b)=(c,b)\) for some \(c\in{\fix(\sigma)}^\ell\).
	Setting \(K:=k(b)\) we have that \(a\ind_kK\) and \(a\in K(c)\), the latter 
	witnessed by \(g_b^{-1}\).
	This shows that \(\qftp(a/k)=p\) is qf-internal to \(\fix(\sigma)\).

	Suppose, now, that \(\characteristic(k)=0\) and~\(p\) is qf-internal to 
	\(\fix(\sigma)\).
	Let \(a\models p\).
	Using condition~(iv) of Proposition~\ref{prop:equivalences-qfint}
	we have \(K\supseteq k\) with \(a\ind_kK\), and an \(\ell\)-tuple~\(c\) 
	from \(\fix(\sigma)\) such that \(K(a)=K(c)\).
	(This is where characteristic zero is being used, we do not have to take 
	the perfect closure.) Moreover, part of condition~(iv) of 
	Proposition~\ref{prop:equivalences-qfint} tells us that we can take \(K\) 
	to be of the form \(K=k(b)\), where \(r:=\qftp(b/k)\) is rational.
	Let \((Z,\psi)\) be the rational \(\sigma\)-variety over~\(k\) whose 
	quantifier-free generic type is~\(r\).
	Let \(Y:=\loc(c,b/k)\).
	It follows that~\(Y\) is \((\id\times\psi)\)-invariant in~\(\AA^\ell\times 
	Z\).
	Note that \(\loc(a,b/k)=V\times Z\) as \(a\ind_kK\).
	Let~\(g:V\times Z\dto Y\) be the birational map such that \(g(-,b)\) 
	witnesses \(K(a)=K(c)\).
	Note that
	\[
	(\id\times\psi)g(a,b)=(c,\psi(b))
	\]
	and also that
	\begin{eqnarray*}
		g^\sigma(\phi\times\psi)(a,b)
		&=&g^\sigma(\phi(a),\psi(b))\\
		&=&g^\sigma(\sigma(a),\sigma(b))\ \ \text{ as \((a,b)\) are 
		\(\sharp\)-points} \\
		&=&\sigma(g(a,b))\\
		&=&\sigma(c,b)\\
		&=&(c,\psi(b)).
	\end{eqnarray*}
	As \((a,b)\) is Zariski generic in \(V\times Z\) over~\(k\), this means 
	that
	\((\id\times\psi)g=g^\sigma(\phi\times\psi)\).
	So~\(g\) is equivariant.
	We have thus produced a trivialisation.
\end{proof}

The above proof gives us a little more that it is worth extracting for later use:

\begin{corollary}\label{cor:ds-iso}
	Suppose \((V,\phi)\) is a rational \(\sigma\)-variety over an algebraically 
	closed difference field~\((k,\sigma)\) of characteristic zero.
	If \((V,\phi)\) is isotrivial then there exists a trivialisation 
	where~\(Z\) is an invariant subvariety of 
	\((V^n\times\AA^m,\phi\times\id)\), for some \(n,m\geq 0\), that projects 
	dominantly onto \(V^n\), and such that~\(\psi\) is the restriction of 
	\(\phi\times\id\) to~\(Z\).
\end{corollary}

\begin{proof}
	Let~\(p\) be the generic quantifier-free type of~\((V,\phi)\) over~\(k\).
	By Proposition~\ref{prop:ds-iso}(a), \(p\) is qf-internal to 
	\(\fix(\sigma)\).
	Now, the proof of Proposition~\ref{prop:ds-iso}(b) constructs a 
	trivialisation of~\((V,\phi)\) that has the additional property we are 
	seeking.
	Indeed, condition~(iv) of Proposition~\ref{prop:equivalences-qfint} ensures 
	that the tuple~\(b\) used in that construction is of the form 
	\(b=(a_1,\dots,a_n,d)\) where~\(a_1,\dots,a_n\) are independent 
	realisations of~\(p\) and~\(d\) is an \(m\)-tuple from \(\fix(\sigma)\).
	It follows that \(Z=\loc(b/k)\) and~\(\psi\) are of the desired form.
\end{proof}

\begin{question}
	The statement of Corollary~\ref{cor:ds-iso}
	does not mention any model theory, but its proof goes via the model-theoretic 
	arguments of \S\ref{subsect:qfint}.
	Is there a purely algebraic dynamics proof of this result?
	Such a proof might very well extend to arbitrary characteristic.
\end{question}

We also obtain a geometric formulation of Proposition~\ref{nonortho-qfint} 
that may be of independent interest:

\begin{corollary}\label{nonortho-qfint-geo}
	Suppose \((V,\phi)\) is a rational \(\sigma\)-variety over an algebraically 
	closed difference field~\((k,\sigma)\).
	Suppose \(\characteristic(k)=0\).
	The following are equivalent:
	\begin{itemize}
		\item[(i)]
			There is a rational \(\sigma\)-variety \((W,\psi)\) over~\(k\) such 
			that \((V,\phi)\times(W,\psi)\) admits an invariant rational function 
			that is not the pullback of a rational function on $W$.
		\item[(ii)]
			There is a dominant equivariant rational map \((V,\phi)\dto(V',\phi')\) 
			over \(k\) with \((V',\phi')\) isotrivial and positive-dimensional.
	\end{itemize}
\end{corollary}

\begin{proof}
	This is a matter of putting together 
	Propositions~\ref{nonortho-qfint},~\ref{prop:ds-wo}, and~\ref{prop:ds-iso}.

	Let \(p\) be the generic quantifier-free type of \((V,\phi)\).
	Condition~(i) is equivalent to \(p\) being nonorthogonal to 
	\(\fix(\sigma)\), by~\ref{prop:ds-wo}(b).
	By~\ref{nonortho-qfint}, this is in turn equivalent to the existence of a 
	rational map \(p\to q\) where \(q\in S_{\qf}(k)\) is positive-dimensional 
	rational and qf-internal to~\(\fix(\sigma)\).
	Such \(p\to q\) corresponds to a dominant equivariant rational map 
	\((V,\phi)\dto(V',\phi')\).
	That~\(q\) is positive-dimensional is equivalent to \(V'\) being 
	positive-dimensional, and that~\(q\) is qf-internal to~\(\fix(\sigma)\) is 
	equivalent to \((V'\phi')\) being isotrivial.
	The latter is by~\ref{prop:ds-iso} as we are in characteristic zero.
\end{proof}

\begin{remark}
	Like Corollary~\ref{cor:ds-iso}, Corollary~\ref{nonortho-qfint-geo} 
	does not mention any model theory, but we have given a model-theoretic proof.
	In this case, however, we do see an algebraic-geometric approach, along the following lines:
	After taking projective closures, a rational function~\(\lambda\) on \(V\times W\) induces a rational map \(f_\lambda\) from~\(V\) to the Hilbert scheme of rational functions on~\(W\), given by \(a\mapsto \lambda(a,-)\), whose image we can take to be~\(V'\).
	If \(\lambda\) is invariant for \(\phi\times\psi\), and assuming that \(\psi\) is birational, we can give \(V'\) a \(\sigma\)-variety structure \(\phi'\) defined by precomposition with \(\psi^{-1}\).
	It is then not hard to verify that \(f_\lambda:(V,\phi)\dto (V',\phi')\) is equivariant and that \((V',\phi')\) is isotrivial. Finally,  if~\(\lambda\) does not arise as the pullback of a rational function on~\(W\) then \(V'\) will be positive-dimensional.	\end{remark}

\bigskip
\section{The binding group theorems}
\label{sect:bg}

\noindent
In this section we prove Theorems~\ref{thm:bgrds} and~\ref{thm:bgrt}.
For that purpose we now restrict entirely to characteristic zero.
We will point out where this is used along the way.

So our assumptions on our base difference field $(k,\sigma)$ are that it is inversive, algebraically closed, and of characteristic zero.
Work in a 
sufficiently saturated model \((\U,\sigma)\models\acfa_0\) extending 
\((k,\sigma)\).
Denote the fixed field of $(\U,\sigma)$ by \(\C\).

Fix also a rational type \(q(x)\in S_{\qf}(k)\) that is qf-internal 
to~\(\C\).

Recall from Definition~\ref{def:bgrt} that the \emph{quantifier-free binding 
group of~\(q\) with respect  to~\(\C\)}, denoted by \(\aut_{\qf}(q/\C)\), is 
the (abstract) subgroup of permutations~\(\alpha\) of \(q(\U)\) satisfying:
\begin{itemize}
	\item[\((\star)\)] For any quantifier-free formula \(\theta(x,y)\) over 
		\(k\), any tuple \(a\) of realisations of~\(q\), and any tuple \(c\) of 
		elements of \(\C\),
		\[\models\theta(a,c)\iff\models\theta(\alpha(a),c).\]
\end{itemize}
It is not hard to see that the set of such permutations does form a subgroup.

\begin{remark}\label{rem:bg}
	\begin{itemize}
		\item[(a)]
			The binding group can be understood as an automorphism group for a 
			certain auxiliary two-sorted structure~\(\Q\), whose sorts are 
			\(q(\U)\) and \(\C\) and where the language is made up of a predicate 
			symbol for each relatively quantifier-free \(k\)-definable subset of 
			\({q(\U)}^n\times\C^m\) in \((\U,\sigma)\).
			Let \(\aut(\Q/\C):=\{\alpha\in\aut(\Q):\alpha|_{\C}=\id_{\C}\}\).
			It can be easily seen that the map \(\aut(\Q/\C)\to \aut_{\qf}(q/\C)\), 
			given by restriction to \(q(\U)\), is an isomorphism of groups that 
			preserves the action on \(q(\U)\).
		\item[(b)]
			As $q$ is rational we need only consider quantifier-free 
			$\lring$-formulas in~\((\star)\).
		\item[(c)]
			The binding group is a birational invariant in the sense that any 
			birational equivalence, \(\gamma:q\to\widehat{q}\), between rational 
			types, lifts canonically to an isomorphism of group actions,
			\[
			\gamma^*:\aut_{\qf}(q/\C)\to\aut_{\qf}(\widehat{q}/\C),
			\]
			given by \(\alpha\mapsto\gamma\alpha\gamma^{-1}\).
	\end{itemize}
\end{remark}

We focus, first of all, on proving that \(\aut_{\qf}(q/\C)\), along with its 
action on \(q(\U)\), has a quantifier-free definable avatar in 
\((\U,\sigma)\).
This is the main clause of Theorem~\ref{thm:bgrt}, and will occupy us for 
most of the section.
Our construction is informed by those of Hrushovski~\cite{udigeneral} and the 
first author~\cite{Moshe},
but with particular attention paid to the birational geometric setting in 
which we find ourselves.

Since~\(q\) is rational it is the generic quantifier-free type of a rational 
\(\sigma\)-variety \((V,\phi)\) over~\(k\).
Since~\(q\) is qf-internal to~\(\C\), and we are in characteristic zero, 
Proposition~\ref{prop:ds-iso} tells us that \((V,\phi)\) is isotrivial.
In fact, by Corollary~\ref{cor:ds-iso}, we have a trivialisation
\[
\xymatrix{
(V\times\widetilde{Z},\phi\times\widetilde\psi)
\ar[dr]\ar@{-->}[rr]^{\widetilde{g}}_{\isom}&& 
(\widetilde{Y},\id\times\widetilde\psi)\ar[dl]\\
&(\widetilde{Z},\widetilde\psi)
}
\]
where
\(\widetilde Y\subseteq \AA^\ell\times\widetilde Z\) is invariant for 
\(\id\times\widetilde \psi\), and \(\widetilde Z\) is an invariant subvariety 
of \((V^n\times\AA^m,\phi\times\id)\), for some \(n,m\geq 0\), that projects 
dominantly onto \(V^n\), and such that~\(\widetilde \psi\) is the restriction 
of \(\phi\times\id\) to~\(Z\).

We make this trivialisation more canonical by applying 
Proposition~\ref{prop:canfam} to the graph of~\(\widetilde g\).
Note that this graph is an invariant subvariety by equivariance -- see 
Lemma~\ref{lem:graph}.
What~\ref{prop:canfam} yields is a dominant equivariant \(\mu:(\widetilde 
Z,\widetilde\psi)\dto(Z,\psi)\) such that for general \(e,e'\in \widetilde 
Z\), \(\mu(e)=\mu(e')\) if and only if \(\widetilde g_e=\widetilde g_{e'}\).
It follows that~\(\widetilde g\) descends to an equivariant birational map
\(g:(V\times Z,\phi\times \psi)\dto (Y,\id\times \psi)\) over \((Z,\psi)\),
where~\(Y\) is now the invariant subvariety of \((\AA^\ell\times 
Z,\id\times\psi)\) obtained as the (Zariski closure of the) image 
of~\(\widetilde Y\) under \(\id\times\mu\).
We thus obtain a trivialisation
\[
\xymatrix{
(V\times Z,\phi\times \psi)\ar[dr]\ar@{-->}[rr]^{ g}_{\isom}&& 
(Y,\id\times \psi)\ar[dl]\\
&(Z,\psi)
}
\]
such that  the family of birational maps
\((g_e:V\dto Y_e:e\in Z)\)
is canonical in the sense that if \(g_e=g_{e'}\) then \(e=e'\), for general 
\(e,e'\in Z\).

We will use both the canonicity of this family of birational maps and the 
fact that \((Z,\psi)\) is the image of \((\widetilde 
Z,\widetilde\psi)\subseteq(V^n\times\AA^m,\phi\times\id)\).

\begin{remark}\label{rem:char0}
	Our dependence on characteristic zero ends here.
	That is, given a trivialisation of \((V,\phi)\) with \((Z,\psi)\) of the 
	above form -- so both canonical and the image of of some \((\widetilde 
	Z,\widetilde\psi)\subseteq(V^n\times\AA^m,\phi\times\id)\) --  the rest of 
	our construction of the binding group, and hence of 
	Theorems~\ref{thm:bgrds} and~\ref{thm:bgrt} go through in any 
	characteristic.
\end{remark}

Let \(r\) be the generic quantifier-free type of \((Z,\psi)\) over~\(k\).
Let us first observe that~\(\aut_{\qf}(q/\C)\) acts on \(r(\U)\) as well:

\begin{lemma}\label{lem:gonr}
	There is an action of~\(\aut_{\qf}(q/\C)\) on \(r(\U)\) such that, for any $e\models r$ and \(\alpha\in\aut_{\qf}(q/\C)\), 
	\(
	g_e(b)=g_{\alpha(e)}(\alpha(b))
	\)
	for all \(b\models q\) with \(b\in\dom{g_e}\subseteq V\).
	Moreover, $\alpha(e)$ is the unique realisation of~$r$ with this property.
\end{lemma}

\begin{proof}
	Because of the the dominant equivariant rational map 
	\(\mu:(\widetilde Z,\widetilde\psi)\dto(Z,\psi)\),
	and the nature of \((\widetilde Z,\widetilde\psi)\),
	realisations of~\(r\) are of the form \(\mu(\overline a,d)\), for some 
	\(\overline a\models q^{(n)}\) and \(d\in\C^m\) such that $(\overline a, d)$ is generic in $(\widetilde Z,\widetilde\psi)$.
	The action we have in mind, for \(\alpha\in\aut_{\qf}(q/\C)\), is 
	\(\mu(\overline a,d)\mapsto \mu(\alpha \overline a,d)\).
	Note that~(\(\star\)) ensures that \((\alpha\overline a,d)\) is again 
	generic in \((\widetilde Z,\widetilde \psi)\), and hence \(\mu(\alpha 
	\overline a,d)\) is again a realisation of~\(q\).
	This is well-defined because~(\(\star\)) also ensures that if 
	\(\mu(\overline a,d)=\mu(\overline a',d')\) then
	\(\mu(\alpha\overline a,d)=\mu(\alpha\overline a',d')\).

	Fix \(e=\mu(\overline a,d)\models r, b\models q\) with \(b\in\dom{g_e}\), 
	and  \(\alpha\in\aut_{\qf}(q/\C)\).
	Then \((b,e)\in{(V\times Z,\phi\times\psi)}^\sharp\) and hence 
	\(g(b,e)=(g_e(b),e)\in{(Y,\id\times\psi)}^\sharp\) as~\(g\) is equivariant.
	It follows that \(g_e(b)=:c\in Y(\C)\).
	Applying~(\(\star\)) to the fact that \(g_{\mu(\overline a,d)}(b)=c\) we 
	deduce that
	\(g_{\mu(\alpha\overline a,d)}(\alpha(b))=c\) as well.
	That is, \(g_{\alpha(e)}(\alpha(b))=c\), as desired.
	
	Finally, the uniqueness of $\alpha(e)$ with this property follows from the canonicity of the family of birational maps  given by~$g$.
\end{proof}

Let \(f:(Y,\id\times \psi)\dto(V\times Z,\phi\times\psi)\) be the inverse 
to~\(g\).
So \(f_e:Y_e\dto V\) is the birational inverse to \(g_e\) given by
\(y\mapsto{}\pi_1(f(y,e))\).
The definable copy of~\(\aut_{\qf}(q/\C)\) that we will eventually construct 
will come from identifying elements of~\(\aut_{\qf}(q/\C)\) with birational 
maps of the form \(f_{e'}\circ g_e:V\dto V\) for certain pairs \((e,e')\) of 
realisations of~\(r\).

\begin{proposition}\label{prop1}
	Suppose \(e,e'\) realise \(r\) with \(\qftp(e/\C)=\qftp(e'/\C)\).
	Then \[f_{e'}\circ g_e:V\dto V\] is a birational map that is defined on all 
	realisations of \(q\), and whose restriction to \(q(\U)\), say 
	\(\alpha=\alpha_{e,e'}\), is an element of \(\aut_{\qf}(q/\C)\).


	Conversely, if \(\beta\in\aut_{\qf}(q/\C)\) and \(e\models r\) then 
	\(\qftp(e/\C)=\qftp(\beta(e)/\C)\) and \(\beta=\alpha_{e,\beta(e)}\).  That 
	is, \(\beta=(f_{\beta(e)}\circ g_e)|_{q(\U)}\).
\end{proposition}

\begin{remark}
	Note that we are not claiming that \(g_e\) is defined on all of \(q(\U)\), 
	just that the composition \(f_{e'}\circ g_e\) is.
\end{remark}

\begin{proof}[Proof of~\ref{prop1}]
	First of all, we need to observe that the composition \(f_{e'}\circ g_e\) 
	makes sense.
	Since \(g_e:V\dto Y_e\) and  \(f_{e'}:Y_{e'}\dto V\) are birational maps, 
	it suffices to show that \(Y_e=Y_{e'}\).
	To that end, observe that, as \(e\in{(Z,\psi)}^\sharp\), 
	\({(Y_e)}^\sigma=Y^\sigma_{\psi(e)}\).
	On the other hand, as~\(Y\) is an invariant subvariety of \((\AA^\ell\times 
	Z,\id\times\psi)\), \(Y_e\subseteq Y^\sigma_{\psi(e)}\).
	By dimension considerations, it follows that \({(Y_e)}^\sigma=Y_e\).
	It follows that \(Y_e\) is defined over \(k(e)\cap\C\).
	As $Y_e$ is absolutely irreducible, and $\C$ is pseudo-algebraically closed, we have that \(Y_{e}(\C)\) is Zariski dense in \(Y_{e}\), and similarly for 
	\(Y_{e'}(\C)\).
	It suffices to show, therefore, that \(Y_{e}(\C)=Y_{e'}(\C)\).
	But this is the case as, for any \(c\in\C^\ell\), the statement that \(c\in 
	Y_{e}\) is part of \(\qftp(e/\C)\), and by assumption  
	\(\qftp(e/\C)=\qftp(e'/\C)\).

	We now proceed by a series of claims.

	\begin{claim}\label{onqe}
		\(f_{e'}\circ g_e\) is defined on all realisations of \(q_e\), the 
		nonforking extension of~\(q\) to \(k\cup\{e\}\).
	\end{claim}
	\begin{claimproof}
		If \(a\models q_e\) then it is Zariski generic in \(V\) over \(k(e)\) and 
		hence outside the indeterminacy locus of \(g_e\).
		Moreover, as \(e\in{(Z,\psi)}^\sharp\) and \(g\) is equivariant, it 
		follows that \(c:=g_e(a)\in Y_e(\C)\).
		Since \(f_e\) is defined at \(c\), so is~\(f_{e'}\).
	\end{claimproof}

	Let us denote by \(\alpha\) the restriction of \(f_{e'}\circ g_e\) to 
	realisations of \(q_e\).

	\begin{claim}\label{star-qe}
		Condition~\((\star)\) holds of \(\alpha\) on realisations of \(q_e\).
		That is, for any quantifier-free formula \(\theta(x,y)\) over \(k\), any 
		tuple \(a\) of realisations of~\(q_e\), and any tuple \(c\) of elements 
		of \(\C\),
		\[\models\theta(a,c)\iff\models\theta(\alpha(a),c).\]
		In particular, if \(a\models q_e\) then \(\alpha(a)\models q\)
	\end{claim}

	\begin{claimproof}
		Taking negations it suffices to prove the left to right direction.
		Let \(d:=g_e(a)\), which, as we saw in the proof of the last claim, is a 
		tuple of elements of \(Y_e(\C)\).
		Then \(f_e(d)=a\).
		So \(\models\theta(a,c)\) tells us that \(\models\theta(f_e(d),c)\), so 
		that \(\models\theta(f_{e'}(d),c)\).
		But \(f_{e'}(d)=f_{e'}(g_e(a))=\alpha(a)\), as desired.
	\end{claimproof}

	\begin{claim}\label{onV}
		Suppose \(u,u'\models r\) such that \(u\ind e\) and 
		\(\qftp(eu/\C)=\qftp(e'u'/\C)\).
		Then \(f_{e'}\circ g_e=f_{u'}\circ g_u\) as birational maps on \(V\).
	\end{claim}

	\begin{claimproof}
		We already know that \(Y_e=Y_{e'}\).
		For the same reasons, \(Y_u=Y_{u'}\).
		It follows that \(g_{u}\circ f_{e}\) and \(g_{u'}\circ f_{e'}\) are both 
		birational maps from \(Y_e\) to \(Y_u\).
		Moreover, they agree on the \(\C\)-points since 
		\(\qftp(eu/\C)=\qftp(e'u'/\C)\).
		But, as $\C$ is pseudo-algebraically closed, the \(\C\)-points are Zariski dense, and we have \(g_{u}\circ 
		f_{e}=g_{u'}\circ f_{e'}\).
		Now, let~\(a\) realise \(q_{eu}\), the nonforking extension of~\(q\) to 
		\(k\cup\{e,u\}\).
		Note that \(g_u\) is defined  at \(a\) because \(a\models q_u\), and
		\begin{eqnarray*}
			g_{u}(a)
			&=&
			g_{u}(f_{e}(g_e(a)))\\
			&=&
			g_{u'}(f_{e'}(g_e(a)))\ \ \text{ as }g_{u}\circ f_{e}=g_{u'}\circ 
			f_{e'}\\
			&=&
			g_{u'}(\alpha(a)).
		\end{eqnarray*}
		Hence
		\[f_{e'}(g_e(a))=\alpha(a)=f_{u'}(g_{u'}(\alpha(a)))=f_{u'}(g_u(a)).\]
		That is, \(f_{e'}\circ g_e\) and \(f_{u'}\circ g_u\) agree on 
		realisations of \(q_{eu}\), and so by Zariski-density on all of \(V\).
	\end{claimproof}

	\begin{claim}\label{onq}
		\(f_{e'}\circ g_e\) is defined on all realisations of \(q\).
		Moreover, if \(\alpha=\alpha_{e,e'}\) now denotes the restriction of 
		\(f_{e'}\circ g_e\) on \(q(\U)\), then \(\alpha\in\aut_{\qf}(q/\C)\).
	\end{claim}

	\begin{claimproof}
		Suppose \(a\models q\).
		Choose \(u\models r\) with \(u\ind ea\), and \(u'\models r\) such that 
		\(\qftp(eu/\C)=\qftp(e'u'/\C)\).
		Then, by Claim~\ref{onqe}, \(f_{u'}\circ g_u\) is defined on the 
		realisations of \(q_u\), and hence on~\(a\), while, by Claim~\ref{onV}, 
		\(f_{e'}\circ g_e=f_{u'}\circ g_u\) as birational maps on~\(V\).
		So \(f_{e'}\circ g_e\) is defined at~\(a\).

		For the moreover clause, suppose \(\bar a\) is a tuple of realisations 
		of~\(q\), \(c\) a tuple of elements of \(\C\),  and \(\theta(x,y)\) a 
		quantifier-free formula over \(k\).
		Now let \(u\models r\) with \(u\ind e\bar a\), and \(u'\models r\) such 
		that \(\qftp(eu/\C)=\qftp(e'u'/\C)\).
		Then, by Claim~\ref{star-qe} applied to \(f_{u'}\circ g_u\) restricted to 
		\(q_u\), we have 
		\(\models\theta(\bar{a},c)\iff\models\theta(f_{u'}(g_u(\bar{a})),c)\).
		Since \(f_{u'}\circ g_u=\alpha\) by Claim~\ref{onV}, this shows that 
		\(\alpha\in\aut_{\qf}(q/\C)\).
	\end{claimproof}

	Claim~\ref{onq} is the main clause of the Proposition.

	Finally, for the converse direction suppose \(\beta\in\aut_{\qf}(q/\C)\) 
	and \(e\models r\).
	Recall, by Lemma~\ref{lem:gonr}, that~\(\aut_{\qf}(q/\C)\) acts on 
	\(r(\U)\), and this is what we mean by \(\beta(e)\).
	Moreover, by property~(\(\star\)) of~\(\aut_{\qf}(q/\C)\), we get that 
	\(\qftp(e/\C)=\qftp(\beta(e)/\C)\).
	Hence \(\alpha_{e,\beta(e)}\), which is the restriction of 
	\(f_{\beta(e)}\circ g_e\) to \(q(\U)\), is an element 
	of~\(\aut_{\qf}(q/\C)\) by what we have just proved.
	We want to show it agrees with~\(\beta\).
	Fix \(a\models q\)
	and compute
	\begin{eqnarray*}
		\alpha_{e,\beta(e)}(a)
		&=&
		f_{\beta(e)}(g_e(a))\\
		&=&
		f_{\beta(e)}(g_{\beta^{-1}\beta(e)}(\beta^{-1}\beta (a)))\\
		&=&
		f_{\beta(e)}(g_{\beta(e)}(\beta(a)))\ \ \ \ \ \text{ by 
		Lemma~\ref{lem:gonr} applied to \(\beta^{-1}\in\aut_{\qf}(q/\C)\)}\\
		&=&
		\beta(a)
	\end{eqnarray*}
	as desired.

	This completes the proof of Proposition~\ref{prop1}.
\end{proof}

Let
\(X:=\{(e,e'): e,e'\models r, \qftp(e/\C)=\qftp(e'/\C)\}\).

Let~\(\Lambda\) be the set of invariant rational function on \((Z,\psi)\), as 
defined in~\(\S\)\ref{subsect:irf}.

\begin{proposition}\label{prop:claim2}
	Fix \(e,e'\models{}r\).
	Then \((e,e')\in X\) if and only if  \(\lambda(e)=\lambda(e')\) for all 
	\(\lambda\in\Lambda\).
	In particular,~\(X\) is quantifier-free-type-definable over~\(k\).
\end{proposition}

\begin{proof}
	Recall that \(\Lambda\) is the set of rational functions~\(\lambda\) on 
	\(Z\) which, when evaluated at some (equivalently any) \(e\models{}r\), 
	lands in the fixed field.
	The desired result is then just Proposition~\ref{prop:isolateC}; namely, 
	the fact that \(\qftp(e/\C)\) is isolated by \(\qftp(e/k,k(e)\cap\C)\), for 
	any~\(e\) realising a rational type.
\end{proof}

Next, let \(E\) to be the equivalence relation on~\(X\) given by
\[(e,e')E(u,u')\iff f_{e'}\circ g_e=f_{u'}\circ g_u \text{ as birational 
transformations of } V.\]
As~\(E\) is relatively definable (even relatively \(\lring\)-definable), we 
have that \(X/E\) is type-definable.
This will be our type-definable copy of~\(\aut_{\qf}(q/\C)\).

The following summarises what we have so far:

\begin{proposition}\label{prop2}
	There is a type-definable group structure,~\(\G\), on \(X/E\), and a 
	type-definable group action of~\(\G\) on \(q(\U)\), such that the 
	groups~\(\G\) and~\(\aut_{\qf}(q/\C)\), along with their actions on 
	\(q(\U)\), are isomorphic.

	More precisely, the association \((e,e')\mapsto \alpha_{e,e'}\) given by 
	Proposition~\ref{prop1} induces a bijection 
	\(\iota:X/E\to\aut_{\qf}(q/\C)\) with the following additional properties:
	\begin{itemize}
		\item[(a)]
			Let \(R_1\subseteq X^3\) be the relatively \(\lring\)-definable ternary 
			relation given by:
			\(\big((e_1,e_1'),(e_2,e_2'),(e_3,e_3')\big)\in R_1\)
			if and only if
			\[(f_{e_1'}\circ g_{e_1})\circ(f_{e_2'}\circ g_{e_2})= f_{e_3'}\circ 
			g_{e_3}\]
			as birational transformations of~\(V\).
			Then~\(R_1/E\) makes \(X/E\) into a group,~\(\G\), such that  
			\(\iota:\G\to\aut_{\qf}(q/\C)\) is an isomorphism of groups.
		\item[(b)]
			Let \(R_2\subseteq X\times{q(\U)}^2\) be the relatively 
			\(\lring\)-definable relation:
			\[\big((e,e'),a,b\big)\in R_2
			\iff
			(f_{e'}\circ g_{e})(a)=b.\]
			Then, modulo~\(E\), the relation \(R_2\) induces a group action 
			of~\(\G\) on \(q(\U)\) such that 
			\((\iota,\id):(\G,q(\U))\to(\aut_{\qf}(q/\C),q(\U))\) is an isomorphism 
			of group actions.
	\end{itemize}
\end{proposition}

\begin{proof}
	Suppose \((e,e'), (u,u')\in X\).
	Then
	\begin{eqnarray*}
		\alpha_{e,e'}=\alpha_{u,u'}
		&\iff&
		f_{e'}\circ g_e|_{q(\U)}=f_{u'}\circ g_u|_{q(\U)}\ \ \ \text{ by 
		construction of \(\alpha\)}\\ &\iff&
		f_{e'}\circ g_e=f_{u'}\circ g_u\text{ on }V\ \ \ \text{ as \(q(\U)\) is 
		Zariski dense in \(V\)}\\
		&\iff&
		(e,e')E(u,u').
	\end{eqnarray*}
	This shows that we have an induced injective map 
	\(\iota:X/E\to\aut_{\qf}(q/\C)\).
	It is surjective as any \(\beta\in\aut_{\qf}(q/\C)\) is of the form 
	\(\alpha_{e,\beta(e)}\), for any \(e\models r\), by the converse direction 
	of Proposition~\ref{prop1}.

	Parts~(a) and~(b) follow rather easily.
	Fix \((e_1,e_1'),(e_2,e_2'),(e_3,e_3')\in X\).
	Then, as \(q(\U)\) is Zariski dense in \(V\), \[
	\big((e_1,e_1'),(e_2,e_2'),(e_3,e_3')\big)\in R_1
	\iff
	(f_{e_1'}\circ g_{e_1})\circ(f_{e_2'}\circ g_{e_2})=
	f_{e_3'}\circ g_{e_3} \text{ on } q(\U).
	\]
	By Proposition~\ref{prop1}, this says that
	\[\big((e_1,e_1'),(e_2,e_2'),(e_3,e_3')\big)\in R_1
	\iff
	\alpha_{(e_1,e_1')}\alpha_{(e_2,e_2')}=\alpha_{(e_3,e_3')} \text{ in } 
	\aut_{\qf}(q/\C).\]
	Since \(\iota:X/E\to\aut_{\qf}(q/\C)\) is induced by 
	\((e,e')\mapsto\alpha_{e,e'}\), this shows that \(R_1/E\) makes \(X/E\) 
	into a group such that~\(\iota\) becomes an isomorphism of groups.

	Fix \((e,e')\in X\) and \(a,b\models q\).
	Then, by Proposition~\ref{prop1}, and construction,
	\[\big((e,e'),a,b\big)\in R_2
	\iff
	\alpha_{e,e'}(a)=b.\]
	This shows that \(R_2\) induces an action of \(\G=X/E\) on \(q(\U)\) that 
	is isomorphic (via~\(\iota\)) to the action of~\(\aut_{\qf}(q/\C)\) on 
	\(q(\U)\).
	This proves part~(b).
\end{proof}

\subsection{Defining \(\G\)}\label{ss:defining}
All that remains of the main clause of Theorem~\ref{thm:bgrt} is to show 
that~\(\G\) is quantifier-free definable.
Note that we do not even know yet that it is quantifier-free-type-definable: 
the quotient of a quantifier-free definable set by an \(\lring\)-definable 
equivalence relation need not be quantifier-free.
However, we will show eventually, in Proposition~\ref{dbg} below,  
that~\(\G\) is actually the set of \(\sharp\)-points of some 
\(\sigma\)-variety structure on an algebraic group.

To that end, we first show how to construct an algebraic group from the 
purely algebraic (canonical) family
\[
\xymatrix{
V\times Z\ar[dr]\ar@{-->}[rr]^g&&Y\ar[dl]\\ 
&Z}
\]
of birational maps on~\(V\).
This part of the construction occurs entirely in~\(\acf\).

Fix a nonempty Zariski open subset \(Z_0\subseteq Z\) such that:
\begin{itemize}
	\item
		\(Z_0\subseteq\pi_2(\dom g)\) so that \(g_e:V\dto Y_e\) is a birational 
		map with inverse \(f_e:Y_e\dto V\), for each \(e\in Z_0\), and
	\item
		if \(e,e'\in Z_0\) and \(g_e=g_{e'}\) then \(e=e'\).
		It follows in this case that \(f_e=f_{e'}\) implies \(e=e'\) as well.
\end{itemize}
Set
\[T:=\{(e,e')\in Z_0\times Z_0:Y_e=Y_{e'}\}.\]
Note that  \(f_{e'}\circ g_e\) is a \(k(e,e')\)-birational transformation 
of~\(V\) for any \((e,e')\in T\).
So~\(E\) extends naturally from~\(X\) to~\(T\).
That is, we now denote by~\(E\) the equivalence relation on~\(T\) given by
\[(e,e')E(u,u')\iff f_{e'}\circ g_e=f_{u'}\circ g_u,\]
set
\[W:=T/E,\]
and denote by
\[\pi:T\to W\]
the quotient map.
We denote by
\[1\in W\]
the element given by \(1:=\pi(u,u)\) for any \(u\in Z_0\).
It corresponds, of course, to the identity birational transformation 
of~\(V\).
We also have
\[
\inv:W\to W
\]
given by \(\inv\pi(u,u')=\pi(u',u)\).

Note that \(Z_0, T, E, W,\pi, 1\), and~\(\inv\) are all \(\lring\)-definable 
over~\(k\).

\begin{remark}\label{rem:uu'}
	If \(\pi(u,u')=\pi(u,u'')\) then \(u'=u''\).
	Indeed, by definition \(f_{u'}\circ g_u=f_{u''}\circ g_u\) on~\(V\), and 
	hence, as \(g_u\) is birational, \(f_{u'}=f_{u''}\), which in turn implies 
	that \(u'=u''\) by the canonicity of the family. If \(\pi(u,u')=w\) we 
	write \(wu\) for \(u'\) (hence \(u=\inv(w)u'\)).  In particular,~\(u\) 
	and~\(u'\) are \(\lring\)-interdefinable over \(k(w)\).
	We note also that since \(w=\pi(u,wu)\), if \(w_1u=w_2u\) for some 
	\(w_1,w_2\in{}W\), then \(w_1=w_2\).
\end{remark}

Consider the subset \(H_0\subseteq W\)\label{ho} made up of those \(w\in W\) 
such that \(wu\) exists for \(u\in{}Z\) generic over \(w\).
That is, to be more precise,
\begin{itemize}
	\item[]
		for any (equivalently some) Zariski generic \(u\in Z\) over \(k(w)\) 
		there is \(u'\in Z_0\) with \((u,u')\in T\) and \(w=\pi(u,u')\).
\end{itemize}
Note that 
\(H_0\) is \(\lring\)-definable; we can quantify over Zariski generic 
\(u\in{}Z\) using definability of types in \(\acf\).
Indeed, if we let \(r_0(u)\) be the Zariski generic type of~\(Z\) (a stationary \(\lring\)-type in \(\acf\) over~\(k\)), and we let \(\phi(u,w)\) be the \(\lring\)-formula saying that \(w\in W\) and there is \(u'\in Z_0\) with \((u,u')\in T\) and \(w=\pi(u,u')\), then \(H_0\) is defined by the \(\phi\)-definition of \(r_0\).

\begin{lemma}\label{lem:wu}
	Suppose \(w\in H_0\) and \(u\in{}Z\) is Zariski generic over \(k(w)\).  
	Then \(wu\) is also  Zariski generic in~\(Z\) over \(k(w)\).
\end{lemma}

\begin{proof}
	This follows immediately from Remark~\ref{rem:uu'}.
\end{proof}

One consequence of this is that we can switch the order of the defining 
condition of \(H_0\), that is: \(w\in H_0\) if and only if
\begin{itemize}
	\item[]
		for any (equivalently some) Zariski generic \(u'\in Z\) over \(k(w)\) 
		there is \(u\in Z_0\) with \((u,u')\in T\) and \(w=\pi(u,u')\).
\end{itemize}
and again \(u\) is completely determined by \(u'\), and is Zariski generic 
over \(k(w)\).

We obtain an \(\lring\)-definable group:

\begin{propdef}\label{g0group}
	\((H_0,1, \cdot,\inv)\) is an \(\lring\)-definable group where we define \(w_1\cdot 
	w_2\) to be the unique \(w_3\in H_0\) with the property that
	\(w_i=\pi(u_i,u_i')\), for \(i=1,2,3\), for some (equivalently any) 
	\((u_i,u_i')\in T\) such that
	\[
	(f_{u_1'}\circ{}g_{u_1})\circ{}(f_{u_2'}\circ{}g_{u_2})=f_{u_3'}\circ{}g_{u_3}
	\]
	on~\(V\).
\end{propdef}

\begin{proof}
Let \(u\in Z\) be Zariski generic over~\(k\).
As~\(1\) is a \(k\)-point we have that~\(u\) is Zariski generic 
over~\(k(1)=k\) and \(1=\pi(u,u)\), witnessing that \(1\in H_0\).

To see that \(H_0\) is preserved by \(\inv\), fix \(w\in H_0\) and~\(u\) 
Zariski generic in~\(Z\) over \(k(w)\).
By Lemma~\ref{lem:wu}, \(wu\) is also Zariski generic over \(k(w)\).
On the other hand, \(k(\inv(w))\subseteq k(w)\) as \(\inv\) is 
\(\acf\)-definable over~\(k\).
So \(\inv(w)=\pi(wu,u)\) witnesses that \(\inv(w)\in H_0\).

Finally, it remains to show that if \(w_1,w_2\in H_0\) then there is 
\(w_3\in{}H_0\) satisfying \(w_1\cdot w_2=w_3\).
(Uniqueness is immediate by the nature of the equivalence relation~\(E\), 
\(\lring\)-definability is clear from the definitions, as is the fact that 
the group axioms are satisfied.)
Let \(u\in Z\) be Zariski generic over \(k(w_1,w_2)\).
We have \(w_1=\pi(u,w_1u)\) and \(w_2=\pi(\inv(w_2)u,u)\).
Now, as \(f_u=g_u^{-1}\),
\[
(f_{w_1u}\circ g_{u})\circ (f_u\circ g_{\inv(w_2)u})=
f_{w_1u}\circ g_{\inv(w_2)u},
\]
and hence
\(w_1\cdot w_2=\pi(\inv(w_2)u,w_1u)=:w_3\).
By the definition of~\(\cdot\), and uniqueness, we get that \(w_3\) is in the 
\(\lring\)-definable closure of \(k(w_1,w_2)\), so in the perfect closure of 
this field.
By~\ref{rem:uu'}, we know that \(\inv(w_2)u\) is Zariski generic in~\(Z\) 
over \(k(w_1,w_2)\), and hence over \(k(w_3)\).
That is, \(w_3=\pi(\inv(w_2)u,w_1u)\) witnesses that \(w_3\in H_0\).
\end{proof}

The \(\lring\)-definable group \(H_0\) was constructed purely out of the 
algebraic family of birational maps on~\(V\).
However, so far, there is no reason why it should be nontrivial.
In fact, in our situation it is.
This is because of the following proposition, which returns to the difference context.

\begin{proposition}\label{ging0}
	\(\G\) is a subgroup of \(H_0\).
\end{proposition}

\begin{proof}
	We only need to show that \(\G\subseteq H_0\), as the group structures were 
	defined identically on both of them -- see Proposition~\ref{prop2}(a).

	Suppose \(w=\pi(e,e')\in \G\) where \((e,e')\in X\).
	By Proposition~\ref{prop1} we know that
	\(\beta:=\alpha_{e,e'}:=(f_{e'}\circ g_e)|_{q(\U)}\)
	is an element of~\(\aut_{\qf}(q/\C)\).
	Fix \(u\models r\) Zariski generic in~\(Z\) over~\(k(w)\).
	The converse direction of Proposition~\ref{prop1} tells us that
	\(\beta=\alpha_{u,\beta(u)}\) as well.
	Hence
	\((f_{e'}\circ g_e)|_{q(\U)}=(f_{\beta(u)}\circ g_u)|_{q(\U)}\),
	which implies that \(f_{e'}\circ g_e=f_{\beta(u)}\circ g_u\)
	on~\(V\).
	That is,
	\(w=\pi(e,e')=\pi(u,\beta(u))\),
	and the latter witnesses that \(w\in H_0\).
\end{proof}

In fact, \(\G\) lands in a much smaller subgroup of \(H_0\).

\begin{definition}
	For each rational function \(\lambda\in{}k(Z)\), let \(H_\lambda\) be the 
	set of those \(w\in{}H_0\) such that \(\lambda(u)=\lambda(wu)\) for some 
	(equivalently any) Zariski generic \(u\in{}Z\) over \(k(w)\).
	For a subset \(A\subseteq{}k(Z)\), we let 
	\(H_A=\bigcap_{\lambda\in{}A}H_\lambda\).
\end{definition}

\begin{proposition}\label{ginhlambda}
	For each \(\lambda\in{}k(Z)\), \(H_\lambda\) is an \(\lring\)-definable 
	subgroup of \(H_0\) over~\(k\).
	Moreover, \(\G\leq{}H_\Lambda\), where \(\Lambda\), recall, is the ring of 
	invariant rational functions on \((Z,\psi)\).
\end{proposition}

\begin{proof}
	To see that \(\G\subseteq{}H_\Lambda\), fix \(w\in \G\) and \(u\models r\) 
	Zariski generic over \(k(w)\).
	We have just seen, in the proof of Proposition~\ref{ging0}, that 
	\(wu=\beta(u)\) for some \(\beta\in\aut_{\qf}(q/\C)\).
	In particular, \(\qftp(u/\C)=\qftp(wu/\C)\) and so 
	\(\lambda(u)=\lambda(wu)\), for all \(\lambda\in \Lambda\).
	This witnesses that \(w\in{}H_\Lambda\).

	To see that \(H_\lambda\) is a subgroup of~\(H_0\), for any 
	\(\lambda\in{}k(Z)\), we just follow the proof of 
	Proposition~\ref{g0group}.
	Namely, given  \(w_1,w_2\in H_\lambda\), let \(u\in Z\) be Zariski generic 
	over \(k(w_1,w_2)\).
	As \(w_1\in H_\lambda\), we have \(w_1=\pi(u,w_1u)\) and 
	\(\lambda(u)=\lambda(w_1u)\).
	Similarly we have \(w_2=\pi(\inv(w_2)u,u)\) and 
	\(\lambda(u)=\lambda(\inv(w_2)u)\).
	Hence \(w_1\cdot w_2=\pi(\inv(w_2)u,w_1u)\), both \(w_1u\) and 
	\(\inv(w_2)u\) are Zariski generic over \(k(w_1,w_2)\), and 
	\(\lambda(w_1u)=\lambda(\inv(w_2)u)\).
	This witnesses that \(w_1\cdot w_2\in H_\lambda\).
	It is also clear that \(H_\lambda\) is preserved by~\(\inv\).
\end{proof}

By the descending chain condition for groups definable in~\(\acf\), 
\(H_\Lambda\) is also \(\lring\)-definable group.
We will show that~\(\G\) is quantifier-free definable in \(\acfa\) by 
endowing \(H_0\) with an \(\lring\)-definable dynamical structure, 
\(\rho_0:H_0\to H_0^\sigma\), and then showing that 
\(\G={(H_0,\rho)}^\sharp\cap H_\Lambda\).

We need two preparatory lemmas that have to do with the transform of the 
situation by~\(\sigma\).  Note that we have 
\(Y^\sigma\subseteq\AA^\ell\times{}Z^\sigma\) a family of subvarieties of 
\(\AA^\ell\), and
\[
\xymatrix{
V^\sigma\times Z^\sigma\ar[dr]\ar@{-->}[rr]^{g^\sigma} 
&&Y^\sigma\ar[dl]\\
&Z^\sigma}
\]
a family of birational maps on \(V^\sigma\), both parameterised by 
\(Z^\sigma\). We also have \(T^\sigma\subseteq Z_0^\sigma\times Z_0^\sigma\) 
and \(\pi^\sigma:T^\sigma\to W^\sigma=T^\sigma/E^\sigma\).

\begin{lemma}\label{toronto-claim2'}
	Suppose \(e\in Z_0\) is in the domain of~\(\psi\).
	Then \(Y_e=Y^\sigma_{\psi(e)}\).
\end{lemma}

\begin{proof}
	We use the fact that~\(Y\) is an invariant subvariety of 
	\((\mathbb{}A^\ell\times{}Z,\psi\times\id)\).
	Since, for any \(x\in\mathbb A^\ell\), we have that 
	\((x,e)\in\dom{\psi\times\id}\), the invariance tells us that
	\begin{eqnarray*}
		x\in Y_e
		&\implies&
		(e,x)\in Y\\
		&\implies&
		(x,\psi(e))\in Y^\sigma\\
		&\implies&
		x\in Y^\sigma_{\psi(e)}.
	\end{eqnarray*}
	That is, \(Y_e\subseteq Y^\sigma_{\psi(e)}\).
	But as \(e\in Z_0\), we have that \(Y_e\) is birationally equivalent 
	to~\(V\) and  \(Y^\sigma_{\psi(e)}\) is birationally equivalent 
	to~\(V^\sigma\), so that these Zariski closed subsets of \(\mathbb A^\ell\) 
	are irreducible and have the same dimension.
	It must therefore be that \(Y_e=Y^\sigma_{\psi(e)}\).
\end{proof}

\begin{lemma}\label{toronto-claim3}
	Suppose  \(u,u'\in Z\) are Zariski generic over~\(k\), and \((u,u')\in T\).
	Then
	\[\phi \circ f_{u'}\circ g_u=f^\sigma_{\psi(u')}\circ 
	g^\sigma_{\psi(u)}\circ\phi\]
	as rational maps \(V\dto V^\sigma\).
\end{lemma}

\begin{proof}
	There are various things to check to even make sense of the statement.
	First of all, as \(u,u'\in Z\) is Zariski generic over~\(k\) we have that 
	\(u,u'\in Z_0\) so that \(g_u,f_{u'}\) are well-defiined birational maps.
	Moreover, \(u,u'\in\dom\psi\) and, as \(\psi: Z\dto Z^\sigma\) is a 
	dominant rational map, we get that \(\psi(u),\psi(u')\in Z_0^\sigma\) so 
	that \(g^\sigma_{\psi(u)},f^\sigma_{\psi(u')}\) are also well-defined 
	birational maps.
	Finally, to compose things, we need to know that 
	\(Y_{\psi(u)}^\sigma=Y_{\psi(u')}^\sigma\).
	This follows by Lemma~\ref{toronto-claim2'} since \(Y_u= Y_{u'}\).

	The identity itself follows readily from the fact that \begin{eqnarray*}
		g^\sigma_{\psi(u)}\circ\phi
		&=&
		g_u,\ \text{ and}\\
		f^\sigma_{\psi(u')}
		&=&
		\phi\circ f_{u'}
	\end{eqnarray*}
	as rational functions on~\(V\) and \(Y_{u'}\), respectively.
	These, in turn, follow from the fact that \(g:(Z\times 
	V,\psi\times\phi)\dto(Y,\psi\times\id)\), and its inverse~\(f\), are 
	equivariant.
\end{proof}

\begin{remark}\label{psi-iso}
	One consequence of the above proof that is worth pointing out is 
	that~\(\psi:Z\dto Z^\sigma\) is necessarily  generically injective (and 
	hence birational in characteristic zero).
	Indeed, we saw that \(g^\sigma_{\psi(u)}\circ\phi = g_u\) for any~\(u\in 
	\dom\psi\), hence, if \(u_1,u_2\in Z_0\) with \(\psi(u_1)=\psi(u_2)\) then 
	\(g_{u_1}=g_{u_2}\), and so, by canonicity, \(u_1=u_2\).
\end{remark}

We now enrich \(H_0\) with dynamics.

\begin{definition}\label{d:rho_0}
	Let \(\rho_0:H_0\to{}H_0^\sigma\) be defined as follows:
	Given \(w=\pi(u,u')\) in \(H_0\), where \(u,u'\) are Zariski generic 
	over~\(k\), set \(\rho_0(w):=\pi^\sigma(\psi(u),\psi(u'))\).
\end{definition}

\begin{proposition}\label{rho_0}
	\(\rho_0\) is a well-defined \(\lring\)-definable group isomorphism 
	over~\(k\).
\end{proposition}

\begin{proof}
	Note, first of all, that such \(u,u'\) exist by definition of~\(H_0\) and 
	Remark~\ref{rem:uu'}.
	Also, as we have already seen, Zariski genericity ensures that
	\(\psi(u),\psi(u')\in Z_0^\sigma\)
	and
	\((\psi(u),\psi(u'))\in T^\sigma\).
	Hence \(\pi^\sigma(\psi(u),\psi(u'))\) makes sense.
	But there are various things to check:

	\begin{enumerate}
		\item\label{rho_0:a} \(\pi^\sigma(\psi(u),\psi(u'))\) depends only 
			on~\(w\) and not on the choice of \(u,u'\). Suppose \(e,e'\in Z\) is 
			another choice of Zariski generic points over~\(k\) with 
			\(\pi(e,e')=w\). Then \(f_{e'}\circ{}g_e=f_{u'}\circ g_u\). Hence
			\begin{eqnarray*}
				f^{\sigma}_{\psi(u')}g^\sigma_{\psi(u)}\phi
				&=&
				\phi f_{u'}g_u\ \ \ \text{ by~Lemma~\ref{toronto-claim3}}\\
				&=&
				\phi f_{e'}g_e\\
				&=&
				f^{\sigma}_{\psi(e')}g^\sigma_{\psi(e)}\phi\ \ \ \text{ 
				by~Lemma~\ref{toronto-claim3} again.}
			\end{eqnarray*}
			As \(\phi\) is dominant, it follows that
			\[
			f^\sigma_{\psi(e')}\circ g^\sigma_{\psi(e)}=
			f^\sigma_{\psi(u')}\circ g^\sigma_{\psi(u)},
			\]
			which says exactly that
			\((\psi(e),\psi(e'))E^\sigma(\psi(u),\psi(u'))\), namely that
			\[
			\pi^\sigma(\psi(e),\psi(e'))=\pi^\sigma(\psi(u),\psi(u')),
			\]
			as desired.

		\item \(\rho_0\) is injective.
			This is just a matter of noticing that each of the steps in the proof 
			of~\eqref{rho_0:a} above are reversible.

		\item \(\rho_0(w)\in H_0^\sigma\).
			The defining condition for \(H_0^\sigma\) is obtained by 
			applying~\(\sigma\) to the defining condition for \(H_0\).
			That is, \(x\in{}W^\sigma\) is in \(H_0^\sigma\) if and only if
			\begin{itemize}
				\item[]
					For any (equivalently some) Zariski generic \(y\in Z^\sigma\) over 
					\(k(x)\) there is \(y'\in Z_0^\sigma\) with \((y,y')\in T^\sigma\) 
					and \(x=\pi(y,y')\).
			\end{itemize}
			Now, we could have chosen~\(u\) to be Zariski generic in~\(Z\) over 
			\(k(w)\). In which case, \(\psi(u)\) is Zariski generic in~\(Z^\sigma\) 
			over \(k(w)\).  At this point, we already know that \(\rho_0\) is an 
			\(\lring\)-definable function over~\(k\), so that \(\rho_0(w)\) is in 
			the perfect closure of \(k(w)\).  Hence \(\psi(u)\) is Zariski generic 
			in~\(Z^\sigma\) over \(k(\rho_0(w))\).  So 
			\(\rho_0(w)=\pi^\sigma(\psi(u),\psi(u'))\) witnesses that
			\(\rho_0(w)\in{}H_0^\sigma\).

		\item \(\rho_0\) is a group homomorphism.
			Here the group structure on  \(H_0^\sigma\) is the one obtained by 
			transforming the group structure on \(H_0\)  by~\(\sigma\).
			Suppose \(w_1\cdot w_2=w_3\) in \(H_0\).  Write \(w_i=\pi(u_i,u_i')\)
			where \(u_i, u_i'\in Z\) are Zariski generic over~\(k\), for 
			\(i=1,2,3\).
			So 
			\((f_{u_1'}\circ{}g_{u_1})\circ(f_{u_2'}\circ{}g_{u_2})=
			f_{u_3'}\circ{}g_{u_3}\) on~\(V\).
			It follows that
			\[
			\phi\circ f_{u_1'}\circ g_{u_1}\circ f_{u_2'}\circ g_{u_2}=\phi\circ 
			f_{u_3'}\circ g_{u_3}.
			\]
			Now, applying Lemma~\ref{toronto-claim3} repeatedly, we deduce that
			\[
			f^\sigma_{\psi(u_1')}\circ g^\sigma_{\psi(u_1)}\circ 
			f^\sigma_{\psi(u_2')}\circ 
			g^\sigma_{\psi(u_2)}\circ\phi=f^\sigma_{\psi(u_3')}\circ 
			g^\sigma_{\psi(u_3)}\circ\phi.
			\]
			As~\(\phi\) is dominant, we get
			\[
			f^\sigma_{\psi(u_1')}\circ g^\sigma_{\psi(u_1)}\circ 
			f^\sigma_{\psi(u_2')}\circ 
			g^\sigma_{\psi(u_2)}=f^\sigma_{\psi(u_3')}\circ g^\sigma_{\psi(u_3)}
			\]
			which says that \(\rho_0(w_1)\cdot\rho_0(w_2)=\rho_0(w_3)\), as 
			desired.

		\item \(\rho_0\) is an isomorphism.
			We have already seen that it is injective. As \(H_0^\sigma\) is an 
			\(\lring\)-definable group of the same Morley rank and degree as 
			\(H_0\), any injective \(\lring\)-definable homomorphism \(H_0\to 
			H_0^\sigma\) is surjective.\qedhere
	\end{enumerate}
\end{proof}

\begin{proposition}\label{prp:Hlaminv}
	For each each \(\lambda\in\Lambda\), the subgroup \(H_\lambda\) of \(H_0\) 
	is \(\rho_0\)-invariant, in the sense that 
	\(\rho_0(H_\lambda)\subseteq{}H_\lambda^\sigma\).
\end{proposition}
\begin{proof}
	Let \(w\in H_\lambda\) and \(u\) Zariski generic in~\(Z\) over \(k(w)\).
	So \(w=\pi(u,wu)\), and hence \(\rho_0(w)=\pi^\sigma(\psi(u),\psi(wu))\) by 
	definition.
	Note that \(\psi(u)\) is Zariski generic in \(Z^\sigma\) over 
	\(k(\rho_0(w))\) since \(\psi\) is dominant and \(\rho_0\) is 
	\(\lring\)-definable.
	Hence, to show that \(\rho_0(w)\in H_\lambda^\sigma\) it suffices to show 
	that \(\lambda^\sigma\psi(u)=\lambda^\sigma\psi(wu)\).
	But \(\lambda=\lambda^\sigma\psi\) as \(\lambda\) is an invariant rational 
	function on~\((Z,\psi)\), and \(\lambda(u)=\lambda(wu)\) as \(w\in 
	H_\lambda\).
\end{proof}

\begin{definition}
	Let \(H:=H_\Lambda\) and \(\rho:=\rho_0|_H:H\to{}H^\sigma\).
\end{definition}

By the equivalence of categories between \(\lring\)-definable groups and 
algebraic groups, we can give \(H\) the structure of a (possibly not 
connected) algebraic group over~\(k\) such that~\(\rho\) is an isomorphism of 
algebraic groups.
That is, \((H,\rho)\) is a \emph{\(\sigma\)-group} in the sense 
of~\cite{kp2007}, except that \((H,\rho)\) isn't technically a 
\(\sigma\)-variety as we have defined it, because~\(H\) is not necessarily 
irreducible.
Nevertheless, much of our terminology about rational \(\sigma\)-varieties 
makes sense and can be used profitably in this setting.

\begin{proposition}\label{dbg}
	\(\G={(H,\rho)}^\sharp:=\{w\in H:\sigma(w)=\rho(w)\}\)
\end{proposition}

\begin{proof}
	We already know from Proposition~\ref{ginhlambda} that \(\G\leq H\).
	To see that \(\G\subseteq{(H_0,\rho_0)}^\sharp\), fix \(w=\pi(e,e')\in \G\) 
	where \(e,e'\in X\).
	In particular, \(e,e'\models r\), so they are Zariski generic over~\(k\) 
	and \(\sigma(e)=\psi(e)\) and \(\sigma(e')=\psi(e')\).
	Hence
	\begin{eqnarray*}
		\rho_0(w)
		&=&
		\pi^\sigma(\psi(e),\psi(e'))\\
		&=&
		\pi^\sigma(\sigma(e),\sigma(e'))\\
		&=&
		\sigma(\pi(e,e'))\\
		&=&
		\sigma(w)
	\end{eqnarray*}
	as desired.

	For the converse, suppose \(w\in{(H,\rho)}^\sharp\).
	Fix \(e\models r\) Zariski generic in~\(Z\) over \(k(w)\).
	Then \(w=\pi(e,we)\).
	Now,
	\begin{eqnarray*}
		\pi^\sigma(\sigma(e),\sigma(we))
		&=&
		\sigma(w)\\
		&=&
		\rho_0(w)\ \ \ \text{ as }w\in{(H_0,\rho_0)}^\sharp\\
		&=&
		\pi^\sigma(\psi(e),\psi(we))\\
		&=&
		\pi^\sigma(\sigma(e),\psi(we))\ \ \ \text{ as }e\in{(Z,\psi)}^\sharp
	\end{eqnarray*}
	Remark~\ref{rem:uu'}, applied to \(\pi^\sigma\), implies that 
	\(\sigma(we)=\psi(we)\).
	Hence \(we\in{(Z,\psi)}^\sharp\) as well, so that \(we\models r\).
	As \(w\in H_\Lambda\) we have that \(\lambda(e)=\lambda(we)\), for all 
	\(\lambda\in \Lambda\).
	Hence, by Proposition~\ref{prop:claim2}, \((e,we)\in X\).
	So \(w\in X/E=\G\).
\end{proof}

In particular, we have now shown that~\(\G\) is a quantifier-free definable 
group. That, together with~\ref{prop2}, completes the proof of the main 
clause of Theorem~\ref{thm:bgrt}.

\medskip
\subsection{Proof of Theorem~\ref{thm:bgrds}}
Let us recall the notation of that theorem.
We denote by \(\bir(V)=\bir_{\U}(V)\) the group of all birational 
transformations of~\(V\) over~\(\U\).
We set \(\VV:=(V,\phi)\) and consider the collection \(\I_\VV\) of all 
irreducible invariant subvarieties of \((V^r\times\AA^s,\phi\times\id)\) 
over~\(k\) that project dominantly onto each copy of~\(V\), as \(r,s\in\NN\) 
vary.
Setting \(\LL=(\AA^1,\id)\), we then consider the subgroup 
\(\bir(\VV/\LL)=\bir_{\U}(\VV/\LL)\) of \(\bir(V)\) made up of those 
birational transformations that preserve each element of~\(\I_{\VV}\).
More precisely, those \(\delta\in\bir(V)\) such that, for each \(X\subseteq 
V^r\times\AA^s\) in~\(\I_{\VV}\), \begin{itemize}
	\item
		\(X\cap(\dom\delta^r\times\AA^s)\) is nonempty, and
	\item
		\(\delta(X)\subseteq X\).
\end{itemize}
Here, \(\delta\) acts diagonally on~\(V^r\) and trivially on~\(\AA^s\).

Theorem~\ref{thm:bgrds} asserts that \(\bir(\VV/\LL)\) is an algebraic group 
of birational transformations over~\(k\).
This is what we want to prove.

Observe that we have already constructed an algebraic group of birational 
transformations of~\(V\), namely~\(H\).
Indeed, there is a rational map
\[
\theta:H\times V\dto V
\]
such that for every \(w\in H\) and \((u,u')\in T\) with \(\pi(u,u')=w\), we 
have the birational transformation
\[
\theta_w:=f_{u'}\circ g_{u}:V\dto V.
\]
By definition of the group structure given in Definition~\ref{g0group} we 
have that
\[
\theta_1=\id_V, \ \text{ and}
\]
\[
\theta_{w_1}\circ\theta_{w_2}=\theta_{w_1\cdot w_2}\ \text{ for all 
}w_1,w_2\in H.
\]
So \(w\mapsto \theta_w\) makes \(H\) a subgroup of \(\bir(V)\), as 
Definition~\ref{def:agbt} requires.

\begin{lemma}\label{lem:ginh}
	\(\bir(\VV/\LL)\leq H\).
	That is, if \(\delta\in \bir(\VV/\LL)\) then there is \(w\in H\) such that 
	\(\delta=\theta_w\).
\end{lemma}

\begin{proof}
	Recall that our canonical trivialisation
	\[
	\xymatrix{
	(V\times Z,\phi\times \psi)\ar[dr]\ar@{-->}[rr]^{ g}_{\isom}&& 
	(Y,\id\times \psi)\ar[dl]\\
	&(Z,\psi)
	}
	\]
	was induced by a trivialisation
	\[
	\xymatrix{
	(V\times\widetilde{Z},\phi\times\widetilde\psi)
	\ar[dr]\ar@{-->}[rr]^{ \widetilde g}_{\isom}&& 
	(\widetilde{Y},\id\times\widetilde\psi)\ar[dl]\\
	&(\widetilde{Z},\widetilde\psi)
	}
	\]
	via a dominant equivariant \(\mu:(\widetilde 
	Z,\widetilde\psi)\to(Z,\psi)\).
	The original trivialisation had the property that \(\widetilde Z\) is an 
	invariant subvariety of \((V^n\times\AA^m,\phi\times\id)\) over~\(k\), for 
	some \(n,m\geq 0\), that projects dominantly onto \(V^n\), and \(\widetilde 
	\psi\) is the restriction of \(\phi\times\id\) on \(V^n\times\AA^m\).
	In particular, \(\widetilde Z\in\I\).

	Let \(\Gamma\) be the graph of~\(\widetilde g\) over~\(\widetilde Z\).
	That is,
	\[\Gamma:=\{(a, z,\widetilde g_z(a)):a\in V, z\in\widetilde Z\}\subseteq 
	V\times\widetilde Z\times\AA^{\ell}\subseteq V^{n+1}\times\AA^{m+\ell}.\]
	Because~\(\widetilde g\) is equivariant and dominant, \(\Gamma\) is 
	invariant for \(\phi\times\widetilde\psi\times\id\) on \(V\times\widetilde 
	Z\times\AA^{\ell}\), by Lemma~\ref{lem:graph}, and hence for 
	\(\phi\times\id\) on \(V^{n+1}\times\AA^{m+\ell}\).
	Moreover, as \(\Gamma\) projects dominantly onto \(V\times\widetilde Z\), 
	it projects dominantly onto each of the \(n+1\) copies of~\(V\).
	Hence, \(\Gamma\in\I_{\VV}\).

	Let \(L\supseteq k\) be a (finitely generated) field extension over 
	which~\(\delta\) is defined.
	Fix \(u\in Z\) Zariski generic over~\(L\), and write \(u=\mu(b)\) where 
	\(b\in\widetilde Z\) is Zariski generic over~\(L\).
	In particular, \(b=(\overline a,d)\) where \(\overline a\in V^n\) is 
	Zariski generic over~\(L\) and \(d\in\AA^m\).
	Since \(\delta\in\bir(\VV/\LL)\) and \(\widetilde Z\in\I_{\VV}\), we have 
	that \(\delta b:=(\delta\overline a,d)\) is again Zariski generic 
	in~\(\widetilde Z\) over~\(L\).
	Fix \(a\in V\) Zariski generic over \(L(b)\).
	Then, as \(\Gamma\in\I_{\VV}\), we have that
	\((\delta a,\delta b,\widetilde g_b(a))\in\Gamma\),
	so that
	\(\widetilde g_{\delta b}(\delta a)=\widetilde g_b(a)\).
	By Zariski genericity of~\(a\), we conclude that
	\(\widetilde g_{\delta b}\circ\delta=\widetilde g_b\).
	Applying~\(\mu\), it follows that
	\(g_{\mu(\delta b)}\circ\delta=g_u\) as rational maps on~\(V\).
	Letting \(u':=\mu(\delta b)\) and \(w=\pi(u,u')\in W\), we have that
	\(\delta=f_{u'}\circ g_u=\theta_w\).

	We claim that \(w\in H_0\).
	That is, we claim that \(w=\pi(e,e')\) where \(e\) (an in fact~\(e'\)) are 
	Zariski generic in~\(Z\) over \(k(w)\).
	(See page~\pageref{ho} to recall the definition of~\(H_0\).)
	To see this let \(e,e'\) realise the same \(\lring\)-type over~\(L\) but 
	independent from \(w\) over~\(L\).
	Since \(\theta_w=\delta\) is over~\(L\), we still have that
	\(f_{e'}\circ g_e=\theta_w\),
	and so \(\pi(e,e')=w\).
	But now, as~\(u\) was chosen Zariski generic over~\(L\), we have that~\(e\) 
	is Zariski generic over~\(k(w)\).

	Finally, we need to show that \(w\in H\).
	That is, we need to show that, for every invariant rational function 
	\(\lambda\in k(Z)\), \(\lambda(e)=\lambda(e')\).
	Note that there is \((\overline a,d)\in\widetilde Z\) generic over~\(L\) 
	such that \(e=\mu(\overline a,d)\) and \(e'=\mu(\delta\overline a,d)\).
	Indeed, this was the case for~\(u,u'\) by construction, and is part of the 
	\(\lring\)-type over~\(L\).
	Pulling back by~\(\mu\), it suffices to show that for every invariant 
	rational function \(\lambda\in k(\widetilde Z)\), \(\lambda(\overline 
	a,d)=\lambda(\delta\overline a,d)\).
	We may assume that \(\lambda\notin k\).
	So \(\lambda:(\widetilde Z,\widetilde\psi)\dto(\AA,\id)\) is dominant and 
	equivariant, and hence its graph
	\[\Gamma(\lambda)\subseteq\widetilde Z\times\AA\subseteq 
	V^n\times\AA^{m+1}\]
	is an element of \(\I_{\VV}\), by Lemma~\ref{lem:graph}.
	It follows, since \(\delta\in\bir(\VV/\LL)\), that
	\[\delta((\overline a,d),\lambda(\overline a,d))
	=((\delta\overline a,d),\lambda(\overline a,d))\in \Gamma(\lambda).\]
	This means that \(\lambda(\overline a,d)=\lambda(\delta\overline a,d)\), as 
	desired.
\end{proof}

\begin{proof}[Proof of Theorem~\ref{thm:bgrds} (conclusion)]
	Let
	\(G:=\{w\in H:\theta_w\in\bir(\VV/\LL)\}\).
	Note that \(G\) is a (possibly not connected) algebraic subgroup of~\(H\) 
	over~\(k\).
	Indeed, the preservation of any fixed \(X\in\I_{\VV}\) is a Zariski closed 
	condition on~\(w\).
	Lemma~\ref{lem:ginh}, identifies \(\bir(\VV/\LL)\) with~\(G\) and thus 
	completes the proof of Theorem~\ref{thm:bgrds}.
\end{proof}

\medskip
\subsection{The rest of Theorem~\ref{thm:bgrt}}
We now make the connection between \(\aut_{\qf}(q/\C)\) and 
\(\bir(\VV/\LL)\), so between \(\G\) and \(G\), as called for by the ``in 
fact" clause of Theorem~\ref{thm:bgrt}.

\begin{lemma}\label{lem:ging}
	\(\G\leq G\).
\end{lemma}

\begin{proof}
	Fix \(w\in \G\).  We need to show that \(\theta_w\) preserves each member 
	of \(\I_{\VV}\).
	Let \(X\subseteq V^r\times\AA^s\) be in \(\I_{\VV}\).
	So, we have an induced rational dynamics
	\[\varphi:=(\phi\times\id)|_X\]
	on~\(X\), such that the first~\(r\) co-ordinate projections 
	\((X,\varphi)\to (V,\phi)\) are dominant equivariant maps.
	Hence, if \(b=(a_1,\dots,a_r,c_1,\dots,c_s)\) is a generic point of 
	\((X,\varphi)\) over~\(k\), then each \(a_i\models q\).
	Since \(w\in G\), we have that \(\theta_w\) restricts to an element of 
	\(\aut_{\qf}(q/\C)\).
	Hence \(\theta_w\) is defined at~\(b\) and the defining 
	condition~(\(\star\)) of the binding group ensures that \(\theta_w(b)\in 
	X\).
	It follows that \((\theta_w\times\id)\) preserves~\(X\), as desired.
\end{proof}

Next we need to show that \(\rho:H\to H^\sigma\) restricts to \(G\to 
G^\sigma\).
This will follow from the following:

\begin{lemma}\label{theta-equivariant}
	\(\theta:(H\times V,\rho\times\phi)\dto(V,\phi)\) is equivariant in the 
	sense that \(\phi\circ\theta=\theta^\sigma\circ(\rho\times\phi)\) as rational maps 
	\(H\times V\dto V^\sigma\).
\end{lemma}

\begin{proof}
For readability we drop the symbol~$\circ$ and denote composition by concatenation.
	It suffices to show that, for each \(w\in H\),
	\[
	\phi\theta_w=\theta^\sigma_{\rho(w)}\phi
	\]
	as rational maps \(V\dto V^\sigma\).
	Letting \(w=\pi(u,u')\) where \(u,u'\in Z\) are Zariski generic over~\(k\), 
	we have that \(\rho(w)=\pi^\sigma(\psi(u),\psi(u'))\) by how~\(\rho\) is 
	defined in Definition~\ref{d:rho_0}.
	So we have that \(\theta_w=f_{u'} g_u\) and 
	\(\theta^\sigma_{\rho(w)}=f^\sigma_{\psi(u')} g^\sigma_{\psi(u)}\)
	and we are trying to prove that
	\[\phi f_{u'} g_u=f^\sigma_{\psi(u')} g^\sigma_{\psi(u)}\phi,\]
	which is exactly Lemma~\ref{toronto-claim3}.
\end{proof}

\begin{proposition}\label{prop:ginheq}
	\(\rho:H\to H^\sigma\) restricts to an isomorphism \(G\to G^\sigma\).
\end{proposition}

\begin{proof}
	Fixing \(w\in G\) we need to show that \(\rho(w)\in G^\sigma\).
	That is, given \(X\subseteq V^r\times\AA^r\) in~\(\I_{\VV}\), we need to 
	show that \(\theta^\sigma_{\rho(w)}\) preserves~\(X^\sigma\).
	Fix \((a,d)\in X\) Zariski generic over~\(k(w)\), where 
	\(a=(a_1,\dots,a_r)\) is tuple of generic points of~\(V\) and\(d\in\AA^s\).
	Since \(\phi\) restricts to a dominant rational map from~\(X\) to 
	\(X^\sigma\),
	we have that \((\phi (a),d)\in X^\sigma\) is Zariski generic over~\(k(w)\).
	And, because \(\theta_{\rho(w)}^\sigma\phi=\phi\theta_w\) by 
	Lemma~\ref{theta-equivariant}, we get
	\[\theta_{\rho(w)}^\sigma(\phi(a),d) = (\phi\theta_w(a),d).\]
	Since \(w\in G\), \(\theta_w\in\bir(\VV/\LL)\), and so \((\theta_w(a),d)\in 
	X\).
	Hence \((\phi\theta_w(a),d)\in X^\sigma\).
	We have shown that \(\theta_{\rho(w)}^\sigma(\phi(a),d)\in X^\sigma\) for a 
	Zariski generic point \((\phi(a),d)\in X^\sigma\) over \(k(w)\).
	This implies that
	\(\theta_{\rho(w)}^\sigma(X^\sigma)\subseteq X^\sigma\), as desired.
\end{proof}

\begin{proof}[Proof of Theorem~\ref{thm:bgrt} (conclusion)]
	It remains only to show that \(\G={(G,\rho|_G)}^\sharp\).
	But Proposition~\ref{dbg} tells us that \(\G={(H,\rho)}^\sharp\) and 
	Lemma~\ref{lem:ging} tell us that \(\G\leq G\).
	From this the result follows.
\end{proof}

\bigskip
\section{Some applications}\label{sect:applications}
\noindent
In this final section we describe some applications of our binding group 
theorems.

We continue to work over a fixed algebraically closed inversive difference field \((k,\sigma)\) of characteristic zero, and a sufficiently saturated  model \((\U,\sigma)\models\acfa_0\) 
extending \((k,\sigma)\), with \(\C=\fix(\sigma)\).
As mentioned in Remark~\ref{rem:char0}, our only use of characteristic zero 
is to deduce that an isotrivial \(\sigma\)-variety has a trivialisation of a 
particularly useful form.

\medskip
\subsection{The autonomous case}\label{subsect:aut}
Here we recover the main results of~\cite{bms} by restricting to the autonomous case when \(k\subset\C\).
So rational $\sigma$-varieties over~$k$ are rational dynamical systems, $(V,\phi)$ where $\phi$ is a dominant rational self-map.

Following~\cite{ch1} and~\cite{bms} we say that a rational dynamical system 
\((V,\phi)\) is \emph{translational} if~\(\phi\) comes from the action of an 
algebraic group; that is, if there is a faithful algebraic group action 
\(\theta:H\times V\to V\) over~\(k\) such that~\(\phi\) agrees with 
\(\theta_h\) for some \(h\in H(k)\).
The following 
recovers~\cite[Corollary A]{bms}
as a corollary of our binding group theorem.

\begin{theorem}\label{thm:translational}
	Suppose \((V,\phi)\) is an isotrivial rational dynamical system.
	Then $(V,\phi)$ is birationally equivalent to a translational dynamical system.
	If~$\phi$ is an automorphism of~$V$, then $(V,\phi)$ itself is translational.\footnote{This theorem, and indeed a general possibly non-autonomous version of it, was announced by Zo\'e Chatzidakis and Ehud Hrushovski in 2008. See~\cite[$\S$1.13]{ch1}.
It is our understanding that a proof is to be included in an upcoming paper.}
\end{theorem}

\begin{proof}
	Let \(q\in S_{\qf}(k)\) be the quantifier-free generic type of 
	\((V,\phi)\).
	It is qf-internal to~\(\C\) by Proposition~\ref{prop:ds-iso}.
	
	We claim that \(\phi|_{q(\U)}\in\aut_{\qf}(q/\C)\).
	Indeed, if \(a\in{(V,\phi)}^\sharp\cap\dom\phi\) then
	\begin{eqnarray*}
		\sigma(\phi(a))
		&=&
		\phi^\sigma(\sigma(a))\\
		&=&
		\phi(\sigma(a))\ \text{ as }\phi^\sigma=\phi\text{ as }k\subseteq\C\\
		&=&
		\phi(\phi(a))\ \text{ as }a\in{(V,\phi)}^\sharp,
	\end{eqnarray*}
	which shows that \(\phi(a)\in{(V,\phi)}^\sharp\).
	Moreover, as \(\phi:V\dto V\) is dominant it takes Zariski generic points 
	to Zariski generic points.
	Hence, if \(a\models q\) then \(\phi(a)\models q\).
	On the other hand, a consequence of isotriviality is that~\(\phi\) is 
	birational.
	Hence, \(\phi|_{q(\U)}\) is at least a permutation of \(q(\U)\).
	To show that it is in \(\aut_{\qf}(q/\C)\) consider a quantifier-free 
	\(\lring\)-formula \(\psi(x,y)\) over \(k\), any tuple \(a\) of 
	realisations of~\(q\), and any tuple \(c\) of elements of \(\C\), and 
	observe that
	\begin{eqnarray*}
		\models\psi(a,c)
		&\iff&
		\models\psi^\sigma(\sigma(a),\sigma(c))\ \text{ as~\(\sigma\) is an 
		\(\lring\)-automorphism}\\
		&\iff&
		\models\psi(\phi(a),c)\ \text{ as~\(a\) is from \({(V,\phi)}^\sharp\) 
		and~\(c\) is from~\(\C\).}
	\end{eqnarray*}
	This proves that \(\phi|_{q(\U)}\in\aut_{\qf}(q/\C)\).
	(As~\(q\) is rational it suffices to verify~(\(\star\)) for 
	\(\lring\)-formulas.)

	We have an algebraic group of birational transformations
	\(\theta:G\times V\dto V\) over~\(k\) given to us by 
	Theorem~\ref{thm:bgrds}.
	Theorem~\ref{thm:bgrt} gives us a (possibly reducible) rational dynamics 
	\(\rho:G\to G\) such that \(\theta\) is equivariant with respect to 
	\(\rho\times \phi\) and~\(\phi\), see Lemma~\ref{theta-equivariant}.
	Moreover, the conclusion of Theorem~\ref{thm:bgrt} is that we can identify 
	\(\aut_{\qf}(q/\C)\) with \({(G,\rho)}^\sharp\) and its action restricted 
	to \(q(\U)\).
	In particular, we have that $\phi|_{q(\U)}=\theta_g|_{q(\U)}$ for some $g\in G(k)$.
	
Let us first consider the case when $\phi$ is an automorphism.
Note that, by definition of an algebraic group of birational transformations of~$V$ (Definition~\ref{def:agbt}), every $h\in G(k)$ determines a birational transformation $\theta_h:V\dto V$, and we will say, for convenience, that $\theta_h$ {\em is an automorphism} if it extends (necessarily uniquely) to an automorphism (namely biregular transformation) of~$V$.
From this it follows that the set of such~$h$ forms a subgroup of $G(k)$.
Note also that $\theta_h$ is an isomorphism if and only if  $\Gamma(\theta_{h})$ -- which recall from~$\S$\ref{sec:ag} is the Zariski closure in $V\times V$ of the set-theoretic graph of $\theta_h$ on its domain -- is the set-theoretic graph of a bijective function from~$V$ to itself.
In particular, the set of such $h$ form the $k$-points of a field-definable subgroup, and hence an algebraic subgroup, say $H\leq G$.
It follows that $\theta|_{H\times V}$ extends (uniquely) to an algebraic group action $\widehat\theta:H\times V\to V$.
Since $\phi|_{q(\U)}=\theta_g|_{q(\U)}$, and $\phi$ is an automorphism, we have that $g\in H(k)$ and $\phi=\widehat\theta_g$.
So $(V,\phi)$ is translational.
	
	Now, let us drop the assumption that $\phi$ is an automorphism, and only try to show that $(V,\phi)$ is birationally equivalent to a translational dynamics.
	We need to upgrade \(\theta\) to an honest regular algebraic group action 
	by automorphisms.
	This can be done because~\(\theta\) makes~\(V\) into a pre-transformation space for~\(G\), in the sense of Zaitsev~\cite{zaitsev}.
	It follows,  by Weil's regularisation theorem (which is~\cite[Theorem~4.9]{zaitsev} in this setting where~$G$ need not be connected) that there is a birational map
	\[
	\gamma:V\dto\widehat{V}
	\]
	and an algebraic group action
	\[
	\widehat\theta:G\times\widehat{V}\to\widehat{V}
	\]
	such that, for each \(w\in G\),
	\[
	\widehat\theta_w=\gamma\theta_w\gamma^{-1}.
	\]
	This is a faithful action: if \(\widehat\theta_w=\id_{\widehat V}\) then 
	\(\theta_w=\id_{V}\) which implies \(w=1\).
	Next, using~\(\gamma\), we can transport the rational dynamics on~\(V\) 
	onto~\(\widehat V\) by setting \[
	\widehat\phi:=
	\gamma\phi\gamma^{-1}:\widehat{V}\dto\widehat{V}.
	\]
	It is \((\widehat V,\widehat\phi)\) that we show is translational.

	We claim, first, that 	
	\(\widehat\theta:(G\times\widehat{V},\rho\times\widehat\phi)\to
	(\widehat{V},\widehat\phi)\) is equivariant.
	Indeed,
	\begin{eqnarray*}
		\widehat\phi\widehat\theta_w
		&=&
		(\gamma\phi\gamma^{-1})(\gamma\theta_w\gamma^{-1})\\
		&=&
		\gamma\phi\theta_w\gamma^{-1}\\
		&=&
		\gamma\theta_{\rho(w)}\phi\gamma^{-1}\ \ \text{ by equivariance of 
		}\theta\\
		&=&
		\gamma\theta_{\rho(w)}{(\gamma)}^{-1}\widehat\phi\\
		&=&
		\widehat\theta_{\rho(w)}\widehat\phi
	\end{eqnarray*}
	for each \(w\in G\), as desired.

	By construction, \(\gamma:(V,\phi)\dto(\widehat V,\widehat\phi)\) is now an 
	equivariant birational map.
	It therefore restricts to a birational equivalence between~\(q\) and the 
	quantifier-free generic type~\(\widehat{q}\) of 
	\((\widehat{V},\widehat\phi)\).
	We obtain an induced isomorphism 
	\[\gamma^*:\aut_{\qf}(q/\C)\to\aut_{\qf}(\widehat{q}/\C),\]
	see Remark~\ref{rem:bg}(c).
	So, \(\widehat\theta\) restricts to an action of \(G:={(G,\rho)}^\sharp\) 
	on \(\widehat{q}(\U)\) that is isomorphic to the action of 
	\(\aut_{\qf}(\widehat{q}/\C)\) on \(\widehat{q}(\U)\).
	Since \(\phi|_{q(\U)}\in\aut_{\qf}(q/\C)\) we also have 
	\(\widehat\phi|_{\widehat{q}(\U)}\in\aut_{\qf}(\widehat{q}/\C)\).
	Hence \(\widehat 
	\phi|_{\widehat{q}(\U)}=\widehat{\theta_h}|_{\widehat{q}(\U)}\) for some 
	\(h\in {(G,\rho)}^\sharp\).
	As~\(\widehat\phi, \widehat \theta, \widehat{q}\) are over~\(k\), it 
	follows that \(h\in H(k)\).
	Finally, as \(\widehat{q}(\U)\) is Zariski dense in~\(\widehat V\), we 
	obtain that  \(\widehat\phi\) agrees with \(\widehat\theta_h\).
	Hence \((\widehat V,\widehat \phi)\) is translational.
\end{proof}

Recall the Zariski dense orbit conjecture:
\emph{if \((V,\phi)\) is a rational dynamical system over~$k$ that has no nonconstant invariant rational functions then there 
exists a \(k\)-point of~\(V\) whose orbit under~\(\phi\) is Zariski dense 
in~\(V\)}.
When~\(k\) is uncountable this is a theorem of Amerik and 
Campana~\cite{amerik-campana}.
For countable~\(k\) it is open in general, though resolved in various cases, 
including when~\(V\) is a smooth projective surface and~\(\phi\) is 
regular~\cite{xie}.
The following special case of the Zariski dense orbit conjecture does not seem to have been addressed in the literature, and follows from Theorem~\ref{thm:translational} by standard arguments:

\begin{corollary}\label{cor:translational}
	Suppose \(\phi:V\to V\) is an automorphism of an algebraic variety 
	over~\(k\) such that \((V,\phi)\) is isotrivial.
	If \((V,\phi)\) admits no nonconstant invariant rational functions then 
	there is \(a\in V(k)\) such that the orbit of~\(a\) under~\(\phi\) is 
	Zariski dense in~\(V\).
\end{corollary}

\begin{proof}
By Theorem~\ref{thm:translational}, $(V,\phi)$ is translational.
For translational algebraic dynamics the truth of the Zariski dense orbit conjecture seems to be well known.
In the absence of an explicit reference, we sketch a proof.

There is an algebraic group action $\theta:G\times V\to V$ such that~$\phi$ agrees with $\theta_g$ for some $g\in G(k)$.
Let~$H$ by the Zariski closure of the subgroup $\langle g\rangle$ generated by~$g$.
As mentioned earlier, because $\U$ is uncountable and algebraically closed, there is $b\in V(\U)$ whose orbit under~$\phi$ is Zariski dense.
It follows that $\O:=H\cdot b$, the $H$-orbit of~$b$, is Zariski dense.
But a Zariski dense orbit of an algebraic group action is always Zariski open.
So $\O$ is Zariski open.
As~$k$ is algebraically closed, $V(k)$ is Zariski dense in~$V$, and hence there exists $a\in V(k)\cap\O$.
Let $Z$ be the Zariski closure of $\{\phi^n(a):n<\omega\}$.
Then $\phi(Z)\subseteq Z$.
As~$\phi$ is an automorphism of~$V$, it follows by noetherianity of the Zariski topology, that $\phi(Z)=Z$.
This in turn implies that~$Z$ is $H$-invariant.
Hence $\O=H\cdot a\subseteq Z$, which forces $Z=V$.
\end{proof}

We also recover the other main result of~\cite{bms}, namely~Corollary~B of that paper, though making use of the following observation occurring
there:

\begin{fact}
\label{fact:translational}
If $(V,\phi)$ is a positive-dimensional translational rational dynamical system over~$k$ then $(V^2,\phi)$ admits a nonconstant invariant rational function.
\end{fact}

\begin{proof}
If $(V,\phi)$ already admits a nonconstant invariant rational function then so does $(V^2,\phi)$.
Otherwise, this fact is precisely~\cite[Proposition~3.2]{bms}.
\end{proof}

\begin{corollary}\label{cor:nonorthbound-auto}
	Suppose \(k\subseteq\C\) is algebraically closed and~\(p\in S_{\qf}(k)\) is 
	a rational type.
	If~\(p\) is nonorthogonal to \(\C\) then \(p^{(2)}\) is not weakly 
	orthogonal to  \(\C\).

	\vspace{\baselineskip}
	\noindent
	{\bf Geometric formulation}:
	Suppose \((V,\phi)\) is a rational dynamical system over~\(k\).
	If  \((V,\phi)\times(W,\psi)\) admits an invariant rational function that 
	is not the pullback of a rational function on \(W\), for some rational 
	\(\sigma\)-variety \((W,\psi)\) over~\(k\), then \((V^{2},\phi)\) admits a 
	nonconstant invariant rational function.
	In particular, if some cartesian power of \((V,\phi)\) admits a nonconstant 
	invariant rational function then already \((V^{2},\phi)\) does.
\end{corollary}

\begin{proof}
	By Proposition~\ref{nonortho-qfint} there is a nonalgebraic rational type 
	\(q\in S_{\qf}(k)\) that is qf-internal to~\(\C\), and a rational map 
	\(p\to q\) over~\(k\).
	And it suffices to show that~\(q^{(2)}\) is not weakly orthogonal 
	to~\(\C\).
	But, by Theorem~\ref{thm:translational}, since we are in the autonomous 
	situation, after possibly replacing~\(q\) with a birationally equivalent 
	quantifier-free type, we may assume that~\(q\) is the generic 
	quantifier-free type of some translational rational dynamical system.
	By Fact~\ref{fact:translational}, translationality implies the cartesian square of this rational dynamical system admits a nonconstant invariant rational function.
	Hence~\(q^{(2)}\) is not weakly orthogonal to~\(\C\), by Proposition~\ref{prop:ds-wo}(a).

	The geometric formulation is obtained by applying the theorem to the 
	generic quantifier-free type~\(p\) of~\((V,\phi)\).
	The assumption on~\((V,\phi)\) tells us that \(p\) is nonorthogonal 
	to~\(\C\), this is Proposition~\ref{prop:ds-wo}(b).
	Hence \(p^{(2)}\) is not weakly orthogonal to~\(\C\), and so~\((V^{2},\phi)\) admits a nonconstant invariant rational function by 
	Proposition~\ref{prop:ds-wo}(a).

	For the ``in particular" clause, suppose \((V^n,\phi)\) admits a 
	nonconstant invariant rational function, and~\(n\) is least such.
	Then \((W,\psi):=(V^{n-1},\phi)\) admits no nonconstant invariant rational 
	functions, and hence those on \((V^n,\phi)=(V,\phi)\times(W,\psi)\) are not 
	pullbacks from \((W,\psi)\), and hence \((V^{2},\phi)\) admits a 
	nonconstant invariant rational function.
\end{proof}

\medskip
\subsection{Dixmier-Moeglin equivalence}
Now we drop the autonomous assumption, and exhibit a new application.
So~\((k,\sigma)\) is an arbitrary algebraically closed difference field of 
characteristic zero.

Given a rational \(\sigma\)-variety \((V,\phi)\) over~\(k\), it is natural to 
ask for conditions that would force there to exist only finitely many maximal 
proper
invariant subvarieties over~\(k\).
A necessary condition is that \((V,\phi\)) admit no nonconstant invariant 
rational function,
as such a rational function would, by taking appropriate level sets, give 
rise to infinitely many codimension~\(1\) invariant subvarieties.
The question of whether this condition is sufficient is sometimes called the 
Dixmier-Moeglin equivalence problem in algebraic dynamics, 
see~\cite[Conjecture~8.5]{BRS} and also~\cite{dme-survey} for a survey of 
Dixmier-Moeglin type problems.
Actually, admitting no nonconstant invariant rational functions is \emph{not} 
sufficient in general (even in the autonomous case), with counterexamples 
given by Henon automorphisms of the affine plane 
(see~\cite[Theorem~8.8]{BRS}).
But we will show that it is sufficient in the (possibly nonautonomous) 
isotrivial case.

But first we need a lemma that is central to how we use binding groups.

\begin{lemma}\label{wotrans}
	Suppose \(q\in S_{\qf}(k)\) is rational and \(q^{(\ell)}\) is weakly 
	orthogonal to~\(\C\).
	Then \(\aut_{\qf}(q/\C)\) acts transitively on \(q^{(\ell)}(\U)\).
\end{lemma}

\begin{proof}
	First of all, Proposition~\ref{prop:isolateC} tells us 
	that~\(\qftp(a/k,k(a)\cap\C)\) isolates \(\qftp(a/\C)\), for any \(a\models 
	q^{(\ell)}\).
	But, as~\(k\) is algebraically closed, weak orthogonality implies that 
	\(k(a)\cap\C\subseteq k\).
	Hence, if \(a_1,a_2\models q^{(\ell)}\) then 
	\(\qftp(a_1/\C)=\qftp(a_2/\C)\).
	Now, recall, from Remark~\ref{rem:bg}(a) that we have the two-sorted 
	auxiliary structure~\(\Q\) whose sorts are \(q(\U)\) and \(\C\) and where 
	the language is made up of a predicate symbol for each relatively 
	quantifier-free \(k\)-definable subset of \({q(\U)}^n\times\C^m\) in 
	\((\U,\sigma)\), for any \(n,m\geq 0\).
	That \(\qftp(a_1/\C)=\qftp(a_2/\C)\) means that, in the structure~\(\Q\), 
	\(\tp_{\Q}(a_1/\C)=\tp_{\Q}(a_2/\C)\).
	Note that \(\Q\) is also sufficiently saturated as it is a reduct of \((\U,\sigma)\) after naming~$k$.
	For the same reason, \(\C\) is also stably embedded 
	in~\(\Q\).
	Hence there is \(\alpha\in\aut(\Q/\C)=\aut_{\qf}(q/\C)\) such that 
	\(\alpha(a_1)=a_2\), as desired.
\end{proof}

\begin{theorem}\label{dme}
	Suppose \(q\in S_{\qf}(k)\) is rational and qf-internal to~\(\C\).
	If~\(q\) is weakly orthogonal to~\(\C\) then it is isolated by a 
	quantifier-free formula.

	\vspace{\baselineskip}
	\noindent
	{\bf Geometric formulation}:
	Suppose \((V,\phi)\) is an isotrivial rational \(\sigma\)-variety 
	over~\(k\) with no nonconstant invariant rational functions.
	Then \((V,\phi)\) has only finitely many maximal proper invariant 
	subvarieties~\(k\).
\end{theorem}

\begin{proof}
	Theorem~\ref{thm:bgrt}, together with Lemma~\ref{wotrans}, gives us a 
	definable group acting relatively definably and transitively on \(q(\U)\).
	Hence \(q(\U)\), being an orbit of this definable group action, is itself a 
	definable set.
	By compactness it is defined by some formula in~\(q\).

	To deduce the geometric formulation we let~\(q\) be the generic 
	quantifier-free type of~\((V,\phi)\).
	Isotriviality implies that~\(q\) is qf-internal to~\(\C\) 
	(Proposition~\ref{prop:ds-iso}), and that there are no nonconstant 
	invariant rational functions implies that~\(q\) is weakly orthgonal 
	to~\(\C\) (Proposition~\ref{prop:ds-wo}(a)).
	So, applying the theorem to~\(q\), we have that \(S:=q(\U)\) is definable.
	Now, note that \(S\) is the complement in~\({(V,\phi)}^\sharp\) of the 
	union of all \({(W,\phi|_W)}^\sharp\) as you range over all proper 
	invariant subvarieties~\(W\) of \((V,\phi)\).
	This is because for any proper subvariety \(W\subset V\), the Zariski 
	closure of \(W\cap{(V,\phi)}^\sharp\) is invariant for \((V,\phi)\).
	It follows from definability of~\(S\) that the above mentioned union is 
	equal to a finite sub-union, which implies that only finitely many of 
	the~\(W\) are maximal.
\end{proof}

\medskip
\subsection{Bounding nonorthogonality}
Our goal in this final subsection is prove a version of 
Corollary~\ref{cor:nonorthbound-auto} above for general (possibly 
nonautonomous) rational \(\sigma\)-varieties.
We therefore drop the assumption that $k\subseteq\C$, and instead work over an arbitrary inversive and algebraically closed difference field $(k,\sigma)$.

It is well known, in stable theories, that a complete type is nonorthogonal 
to a definable set if and only if some Morley power of it is not weakly 
orthogonal.
The version for rational types nonorthogonal to the fixed field in \(\acfa\) 
appeared as Corollary~\ref{cor:o-wo} above.
The question of how high a Morley power one must take was raised in the 
eighties (see~\cite{hrushovski1989almost}) and has been addressed in various 
settings recently, especially for differential-algebraic geometry, 
see~\cite{nmdeg,abred,pperm}.
The main tool in these recent works has been the binding group action.
Now that we have an appropriate quantifier-free binding group theorem we 
obtain the same bound for rational types in \(\acfa_0\).

\begin{theorem}\label{thm:nonorthbound}
	Suppose~\(p\in S_{\qf}(k)\) is rational.
	If~\(p\) is nonorthogonal to \(\C\) then \(p^{(n)}\) is not weakly 
	orthogonal to  \(\C\) where \(n=\dim(p)+3\).

	\vspace{\baselineskip}
	\noindent
	{\bf Geometric formulation}:
	Suppose \((V,\phi)\) is a rational \(\sigma\)-variety over~\(k\) such that,
	for some rational \(\sigma\)-variety \((W,\psi)\) over~\(k\), the cartesian 
	product \((V,\phi)\times(W,\psi)\) admits an invariant rational function 
	that is not the pullback of a rational function on~\(W\).  Then 
	\((V^{n},\phi)\) admits a nonconstant invariant rational function for 
	\(n=\dim V+3\).
	In particular, if some cartesian power of \((V,\phi)\) admits a nonconstant 
	invariant rational function then already \((V^{\dim V+3},\phi)\) does.
\end{theorem}

\begin{proof}
	By Proposition~\ref{nonortho-qfint} there is a nonalgebraic rational type 
	\(q\in S_{\qf}(k)\) that is qf-internal to~\(\C\), and a rational map 
	\(p\to q\) over~\(k\).
	Let \(m:=\dim(q)\leq \dim(p)\).
	It suffices to show that \(q^{(m+3)}\) is not weakly orthogonal to 
	\(\C:=\fix(\sigma)\).
	We assume that \(q^{(m+3)}\) is weakly orthogonal to~\(\C\), and seek a 
	contradiction.

Let $(V,\phi)$ be the isotrivial rational $\sigma$-variety that~$q$ is the generic type of.
Theorem~\ref{thm:bgrds} gives us an algebraic group of birational transformations
	\(\theta:G\times V\dto V\) over~\(k\).
	Theorem~\ref{thm:bgrt} gives us an isomorphism \(\rho:G\to G^\sigma\) making~\(\theta\) equivariant (see Lemma~\ref{theta-equivariant}), and such that the action of $\G:=(G,\rho)^\sharp$ on $q(\U)$ identifies with that of~\(\aut_{\qf}(q/\C)\).
	As in the proof of 
	Theorem~\ref{thm:translational}, the Weil regularisation theorem upgrades~$\theta$ to a regular action by automorphisms, at the expense of replacing~$q$ by a birationally equivalent type.
	
	By Lemma~\ref{wotrans} the diagonal action of~\(\G\) on 
	\(q^{(m+3)}(\U)\) is transitive.
	It follows that the diagonal action of~\(G\) on~\(V^{m+3}\) is generically 
	transitive, in the sense that it has a Zariski dense orbit.
	That is, $G$ acts {\em generically $(m+3)$-transitively} on the $m$-dimensional variety~$V$.
	But this is ruled out by the truth of the Borovik-Cherlin 
	conjecture in \(\acf_0\) from~\cite{nmdeg}.
	Indeed, that $G$ acts generically $(m+2)$-transitively on~$V$ already implies, by~\cite[Theorem~6.3]{nmdeg}, that $(G,V)$ is definably isomorphic to the natural action of~$\operatorname{PGL}_{m+1}$ on~$\PP_m$, and this latter action is not generically $(m+3)$-transitive.
	This is the sought for contradiction.

	The geometric formulation follows exactly as in 
	Corollary~\ref{cor:nonorthbound-auto}, by applying the theorem to the 
	generic quantifier-free type of~\((V,\phi)\) and using 
	Proposition~\ref{prop:ds-wo}.
	The ``in particular" clause also follows exactly as in 
	Corollary~\ref{cor:nonorthbound-auto}.
\end{proof}

In the differential case, where the above bound of \(\dim(p)+3\) was 
established in~\cite{nmdeg} and~\cite[\(\S\)5]{abred}, it is known to be 
sharp.
We ask for the same here:

\begin{question}
	Is the bound in Theorem~\ref{thm:nonorthbound} sharp?
\end{question}


\end{document}